\documentclass{IEEEtran}

\newcommand{\hal}[1]{{}}

\newcommand{\Id}{\mathrm{Id}}
\newcommand{\Rinfty}{\mathscr{R}_\infty}

\usepackage{amssymb,amsmath,amsthm} %,amsbsy,,amsfonts,amscd,mathrsfs}
\usepackage{graphicx,color,mathrsfs,mathabx}

\newcommand{\ffoot}[1]{}

\newcommand{\beq}{\begin{equation}}
\newcommand{\eeq}{\end{equation}}

\newcommand{\Mexp}{\Phi}
\newcommand{\hR}{\Mexp(\mathcal{F})}
\newcommand{\hS}{\mathcal{F}_\infty}

\usepackage[english]{babel}

\newtheorem{Theorem}{Theorem}
\newtheorem{ml}[Theorem]{Lemma}
\newtheorem{example}[Theorem]{Example}
\newtheorem{remark}[Theorem]{Remark}

\newtheorem{corollary}[Theorem]{Corollary}
\newtheorem{Definition}[Theorem]{Definition}
\newtheorem{mpr}[Theorem]{Proposition}
\newtheorem{proposition}[Theorem]{Proposition}

\newcommand{\A}{\mathscr{A}}

\newcommand{\bt}{\begin{Theorem}}
\newcommand{\et}{\end{Theorem}}
\newcommand{\bl}{\begin{lemma}}
\newcommand{\el}{\end{lemma}}
\newcommand{\bp}{\begin{proposition}}
\newcommand{\ep}{\end{proposition}}
\newcommand{\bc}{\begin{corollary}}
\newcommand{\ec}{\end{corollary}}
\newcommand{\bdeff}{\begin{definition}}
\newcommand{\edeff}{\end{definition}}
\newcommand{\brem}{\begin{remark}}
\newcommand{\erem}{\end{remark}}
\newcommand{\bproof}{\begin{IEEEproof}}
\newcommand{\eproof}{\end{IEEEproof}}
\newcommand{\bi}{\begin{itemize}}

\newcommand{\ei}{\end{itemize}}
\newcommand{\bd}{\begin{description}}
\newcommand{\ed}{\end{description}}
\newcommand{\be}{\begin{enumerate}}
\newcommand{\ee}{\end{enumerate}}

\newcommand{\bqn}{\begin{eqnarray}}
\newcommand{\eqn}{\end{eqnarray}}
\newcommand{\eqnn}{\nonumber\end{eqnarray}}
\newcommand{\eqnl}[1]{\label{#1}\end{eqnarray}}

\newcommand{\ba}[1]{\begin{array}{#1}}
\newcommand{\ea}{\end{array}}

\newcommand{\R}{\mathbb{R}}

\newcommand{\N}{\mathbb{N}}

\begin{document}

\title{
A characterization of switched linear control systems with finite $L_2$-gain
}

\author{Yacine~Chitour,
        Paolo~Mason,
        and~Mario~Sigalotti% <-this % stops a space
\thanks{Y. Chitour and P. Mason are with
Laboratoire des Signaux et Syst\`emes (L2S), CNRS - CentraleSupelec - Universit\'e Paris-Sud, 3, rue Joliot Curie, 91192, Gif-sur-Yvette, France,
{\tt yacine.chitour@lss.supelec.fr, paolo.mason@lss.supelec.fr}}
%Mario Sigalotti
\thanks{M. Sigalotti is with INRIA Saclay, Team GECO \&
CMAP, \'Ecole Polytechnique, Palaiseau,   France,
 {\tt mario.sigalotti@inria.fr}}
 \thanks{This research was partially supported by the iCODE Institute, research project of the IDEX Paris-Saclay, and by the Hadamard Mathematics LabEx (LMH) through the grant number ANR-11-LABX-0056-LMH in the ``Programme des Investissements d'Avenir''.}
 }

\maketitle

\begin{abstract}
Motivated by an open problem posed by J.P. Hespanha, we extend the notion of Barabanov norm and extremal trajectory to classes of switching signals that are not closed under concatenation. We use these tools to 
prove that  the finiteness of the $L_2$-gain  is equivalent, for a large set of switched linear control systems, to the condition that the generalized spectral radius associated with any minimal realization of the original switched system is smaller than one.
\end{abstract}


%\IEEEpeerreviewmaketitle

\section{Introduction}
Let $n,m,p$ be positive integers and $\tau$ be a positive real number. Consider the switched linear control system 
\begin{equation}\label{SSSSS}
\dot x=A_\sigma x+B_\sigma u,\quad y=C_\sigma x+D_\sigma u,
\end{equation}
where $x\in \mathbb{R}^n,\ u\in \mathbb{R}^m,\ y\in \mathbb{R}^p$, $A_\sigma, B_\sigma, C_\sigma, D_\sigma$ are matrices of appropriate dimensions and $\sigma$ is in the class $\Sigma_{\tau}$ of piecewise constant signals with dwell time $\tau$ taking values in a fixed finite set $P$ of indices. 
Define the $L_2$-gain as
$$\gamma_2(\tau)\!=\!{\sup}\!\left\{\frac {\|y_{u,\sigma}\|_2}{\|u\|_2}\mid u\in L_2([0,\infty),\mathbb{R}^m)\setminus\{0\},\;\sigma\in\Sigma_{\tau}\!\right\}\!,$$
where $y_{u,\sigma}$ is the output corresponding to the trajectory of the system associated with $u$ and $\sigma$ starting at the origin at time $t=0$. 
%In \cite[Problem 4.1]{Hespanha-unsolved}, J.P. Hespanha asked the following questions: (i) under which conditions is the function $\tau\mapsto \gamma_2(\tau)$ bounded over $(0,\infty)$? (ii) when $\gamma_2$ is not a bounded function over $(0,\infty)$, how to compute $\tau_{\rm min}$, the infimum of the dwell-times $\tau>0$ for which $\gamma_2(\tau)$ is finite? (iii) how regular is $\gamma_2$? 

In recent years the study of the $L_2$-gain for switched linear control systems with minimum dwell time has attracted a significant interest, especially from a computational point of view. The research has been  mainly focused  on the estimation and on the actual computation of the $L_2$-gain. %, neglecting the questions posed in~\cite[Problem 4.1]{Hespanha-unsolved}.

This is a challenging problem, as the $L_2$-gain of a switched linear control system is not just a function of the $L_2$-gain of { modes},
%the single subsystems,
 not even for an arbitrarily large dwell time. It is well known, indeed, 
%(see e.g.~\cite{Hespanha-IEEE}) 
that in general $\lim_{\tau\to \infty} \gamma_2(\tau)>\max_{p\in P}\gamma_2^p$, where $\gamma_2^p$ denotes the $L_2$-gain of the time-invariant control system where $\sigma(\cdot)\equiv p$. One of the first references dealing with this  problem is~\cite{Hespanha-IEEE}, where an algorithm for the computation of the $L_2$-gain in the case of a single switching (or, equivalently, for the computation of $\lim_{\tau\to \infty} \gamma_2(\tau)$) is illustrated. A generalization of this method is introduced in~\cite{GeromelColaneri}, where  it is shown that $\gamma_2(\tau)<\gamma$ if certain linear matrix inequalities, depending on the parameters $\tau$ and $\gamma$, admit a solution.
A further characterization of the  $L_2$-gain (on bounded time intervals) is given in~\cite{Margaliot-Hespanha} in terms of a variational principle. In  the one-dimensional case, this result allows { one} to compute $\gamma_2(\tau)$ by a simple bisection algorithm. 
Note that in all the %previous 
 above-mentioned
results it is assumed that each operating mode of the switched control system satisfies a minimal realization assumption. Moreover the corresponding numerical methods are all related to the %problem of checking the existence of solutions 
study of suitable Riccati equations or inequalities.

The aim of this paper is rather different with respect to the above-cited works. 
The motivation comes from \cite[Problem 4.1]{Hespanha-unsolved}, where J.P. Hespanha asks the following questions: (i) under which conditions is the function $\tau\mapsto \gamma_2(\tau)$ bounded over $(0,\infty)$? (ii) when $\gamma_2$ is not a bounded function over $(0,\infty)$, how to compute $\tau_{\rm min}$, the infimum of the dwell-times $\tau>0$ for which $\gamma_2(\tau)$ is finite? (iii) how regular is $\gamma_2$? 

To address these questions, a rather natural idea is to reduce the issue of verifying the boundedness of the $L_2$-gain for~\eqref{SSSSS} to the stability problem for  the related uncontrolled switched system $\dot x=A_\sigma x$ (without providing any explicit estimate of the $L_2$-gain). Indeed, in the unswitched case, the finiteness of the $L_2$-gain for $\dot x=Ax+Bu,\ y=Cx$ is equivalent of the stability of $\dot x =Ax$ once $(A,B,C)$ is under minimal realization. In the switched framework, it has been suggested, without being formalized (see e.g., \cite{GeromelColaneri,Margaliot-Hespanha}), that  this equivalence still holds true at least whenever each mode of the switched control system satisfies a minimal realization assumption. Note that in~\cite{Chesi} (see also~\cite{GerCol2}) a nonconservative numerical method 
is
%was 
obtained in order to check the stability of switched systems with minimum dwell time.
Therefore, if the specific questions posed by Hespanha %in \cite[Problem~4.1]{Hespanha-unsolved} about $L_2$-gain issues 
could be reduced to a stability problem for a switched system with minimum dwell time, then the  algorithm proposed in \cite{Chesi} could be directly applied in order to compute   the value $\tau_{\min}$.

The problem of checking the finiteness of the $L_2$-gain and the link between this problem and the stability of a corresponding switched system are of interest also for other classes of switching signals than $\Sigma_{\tau}$, for instance signals with average dwell time constraints~\cite{HirataHespanha} or signals satisfying a condition of persistent excitation~\cite{Chitour-S,Sr-Ak}. For this reason we develop in this paper an abstract framework by introducing axiomatic requirements on the family of signals. These assumptions are satisfied by many meaningful families of switching signals, as listed in Section~\ref{s:classes}. % of switching signals. 
Our main result, Theorem~\ref{th1}, states that, assuming without loss of generality that the switched system is in minimal realization, the finiteness of the $L_2$-gain is equivalent to the global asymptotic stability of the uncontrolled switched system  under a suitable uniform observability assumption.  An important aspect of our result is that such a uniform observability assumption is not merely technical, since the equivalence fails in general, as illustrated in Example~\ref{prop-partial}. Note that  the uniform observability assumption is less restrictive compared with the hypothesis that each operating mode is in minimal realization.

One of the key ingredients in the proof of our main result is the use of switching laws expressing the ``most unstable'' behavior  of the uncontrolled switched system. In the case of arbitrary switching  the most unstable behavior of the system is well represented by the notion of extremal trajectory (see e.g.~\cite{barabba88}). For the families $\Sigma_{\tau}$, or for other families of switching signals considered in Section~\ref{s:classes}, this notion cannot be trivially generalized, essentially due to the fact that these classes are not closed  under concatenation. To bypass this issue, our technique consists in  identifying a subset of the switching signals that is large enough to encompass the asymptotic properties of the original uncontrolled switched linear system, and well-behaved with respect to concatenation. 
The associated flows  define a semigroup of matrices,  whose analysis allows us to  describe the asymptotic most unstable behavior of the original switched system. 
Note that for the class of switching signals with minimum dwell time and that of uniformly Lipschitz signals a somehow similar characterization of the most unstable behavior for switched systems may be found in~\cite{Wirth-SIAM}, although in that paper the approach is much more intricate and less flexible, as our assumptions include a larger number of relevant classes of switching signals. 

The paper is organized as follows. In Section~\ref{s:pre} we introduce the notations used in the paper and we list some relevant classes of switching laws that are included in our framework. In Section~\ref{s:barabba} we prove the existence of a ``most unstable behaviour'' by developing the notion of \emph{quasi-Barabanov semigroup}. Section~\ref{s:L2} contains the main result of the paper, Theorem~\ref{th1} and the proof, through Example~\ref{prop-partial}, that the uniform observability assumption cannot just be removed. 
%withdrawn.  
Finally, in Section~\ref{s:rc}, we show a result   establishing the right-continuity of the map $\tau\mapsto \gamma_2(\tau)$. This provides a partial answer to Question (i) and to Question (iii) in \cite[Problem~4.1]{Hespanha-unsolved}.

\section{Preliminaries}
\label{s:pre}
\subsection{Notations}
If $n,m$ are positive integers, the set of $n\times m$ matrices with real entries is denoted 
$M_{n,m}(\mathbb{R})$ and  simply $M_n(\mathbb{R})$ if $n=m$. We use $\Id_n$ to denote the $n\times n$ identity matrix. A norm on $\R^n$ is denoted $\Vert\cdot\Vert$ and similarly for the induced operator norm on $M_n(\mathbb{R})$. 
A subset $\mathscr{M}$ of $M_n(\mathbb{R})$ is said to be irreducible if the only 
subspaces which are invariant for each element of $\mathscr{M}$ are $\{0\}$ or $\mathbb{R}^n$. 
For every $s,t\geq 0$ and $A\in L_\infty([s,s+t],M_n(\R))$, denote by $\Phi_{A}({s+t},s) \in M_n(\R)$ the flow
(or fundamental matrix) 
of 
$\dot x(\tau)=A(\tau)x(\tau)$ from time $s$ to time $s+t$. Given two signals $A_j:[0,t_j]\to \mathscr{M}$, $j=1,2$, we denote by $A_1\ast A_2:[0,t_1+t_2]\to \mathscr{M}$ the \emph{concatenation of $A_1(\cdot)$ and $A_2(\cdot)$}, i.e., the signal coinciding with $A_1(\cdot)$ on $[0,t_1]$ and with 
 $A_2(\cdot\,-t_1)$ on $(t_1,t_1+t_2]$. Similarly, if $\mathcal{A}$ and  $\mathcal{B}$ are two subsets of signals, we use $\mathcal{A}\ast \mathcal{B}$ to denote the set of signals obtained by concatenation of a signal of $\mathcal{A}$ and a signal of $\mathcal{B}$.

 Let $n,p$ and $m$ be positive integers. Consider a switched linear  control system of the type
 \begin{equation}
% \begin{split}%\label{LSCS}
 \dot x(t)=A(t)x(t)+B(t)u(t),\ y(t)=C(t)x(t)+D(t)u(t),\label{LSCS}
%\end{split} 
\end{equation}
 where $x\in \mathbb{R}^n$, $u\in \mathbb{R}^m$, $y\in \mathbb{R}^p$, 
 $(A,B,C,D)$ belongs to a class $\mathcal{T}$ of measurable switching laws taking values in a bounded set of quadruples of matrices $\mathscr{M}\subset M_n(\mathbb{R})\times  M_{n, m}(\mathbb{R})\times  M_{p,n}(\mathbb{R})\times  M_{p,m}(\mathbb{R})$.

For $t\geq 0$ and a switching law $(A,B,C,D)\in \mathcal{T}$, the controllability and observability Gramians in time $t$ are defined respectively as
\begin{eqnarray*}
\int_0^t\Phi_A(0,s)B(s)B(s)^T\Phi_A(0,s)^Tds,\nonumber\\
\int_0^t\Phi_A(s,0)^TC(s)^TC(s)\Phi_A(s,0)ds.%\label{eq:CG-OC}
\end{eqnarray*}

For $T>0$, let $L_2(0,T)$ be the Hilbert space of measurable functions $u:[0,T)\rightarrow\mathbb{R}^m$
with finite $L_2$-norm, i.e., such that $\Vert u\Vert_{2,T}:=\Big(\int_0^T\Vert u(t)\Vert^2dt\Big)^{1/2}$ is finite. 
If $T=\infty$, we simply use $L_2$ and $\Vert u\Vert_2$ to denote respectively the corresponding Hilbert space and $L_2$-norm.  

For $u\in L_2$ and $\sigma=(A,B,C,D)\in \mathcal{T} $, let 
$y_{u,\sigma}$ be the corresponding output of Eq.~\eqref{LSCS} and 
define the $L_2$-gain associated with  $\mathcal{T} $ by
%\begin{equation}\label{L2G}
\[
\gamma_2(\mathcal{T} ):=\sup\left\{
\frac{\Vert y_{u,\sigma}\Vert_2}{\Vert u\Vert_2}\mid u\in L_2\setminus\{0\},\ \sigma\in \mathcal{T}
\right\}.
\]
%\end{equation}
In this paper, we investigate qualitative properties of $\gamma_2(\mathcal{T} )$ and in particular we are interested in finding conditions ensuring its finiteness. We will therefore assume from now on with no loss of generality  that $D(\cdot)\equiv 0$. 

\subsection{Classes of switching functions}\label{s:classes}
We introduce in this section several classes of switching signals contained in $L_\infty([0,\infty),\mathscr{M})$, for some subset $\mathscr{M}$ of a finite-dimensional vector space.
\begin{itemize}
\item $\mathcal{S}^{\rm arb}(\mathscr{M})$ is the class of \emph{arbitrarily switching} signals, i.e., $\mathcal{S}^{\rm arb}(\mathscr{M})=L_\infty([0,\infty),\mathscr{M})$;
\item  $\mathcal{S}^{\rm pc}(\mathscr{M})$ is the class of piecewise constant signals (i.e., signals whose restriction to every finite time-interval admits a finite number of discontinuities and takes finitely many values);
\item  $\mathcal{S}^{{\rm d},\tau}(\mathscr{M})$ is the class of piecewise constant signals {with} \emph{dwell-time} $\tau> 0$, i.e., {such that} the distance between two switching times is at least $\tau$ (notice that $\mathcal{S}^{\rm pc}(\mathscr{M})$ can be identified with  $\mathcal{S}^{{\rm d},0}(\mathscr{M})$); 
\item $\mathcal{S}^{{\rm av-d},\tau,N_0}(\mathscr{M})$ is the class of piecewise constant signals which satisfy the  \emph{$(\tau,N_0)$ average dwell-time} condition with $\tau> 0$ and $N_0$ a positive integer: for every 
 $s,t\geq 0$, the number of switching times in $[s,s+t]$ is bounded from above by $N_0+t/\tau$;
\item for $\mathscr{M}$ of the type $\{M_{{\delta}}=(1-\delta)M_0+{\delta} M_1\mid {\delta}\in [0,1]\}$, 
\begin{align*}
\mathcal{S}^{{\rm pe},T,\mu}(\mathscr{M})=\bigg\{
M_{\alpha} & \mid \alpha\in L_\infty([0,\infty),[0,1]),\\ 
& \int_t^{t+T}\alpha(s)ds\ge\mu\;\forall t\ge0
\bigg\}
\end{align*}
is the class of \emph{$(T,\mu)$-persistently exciting} signals with $0<\mu \leq T$;
\item 
$\mathcal{S}^{{\rm lip},L}(\mathscr{M})$ is the class of Lipschitz signals with Lipschitz constant $L>0$;
\item 
$\mathcal{S}^{{\rm BV},T,\nu}(\mathscr{M})$ is the class of \emph{$(T,\nu)$-BV signals}, i.e., {the signals whose} restriction to every interval of length $T$  
has total variation at most $\nu$, {that is, $M\in \mathcal{S}^{{\rm BV},T,\nu}(\mathscr{M})$ if and only if \[\sup_{\stackrel{t\geq 0,k\in \mathbb{N}}{t=t_0\leq t_1\leq\dots \leq t_k=t+T}}\sum_{i=1}^k \|M(t_i)-M(t_{i-1})\|\leq \nu.\]}
\end{itemize}
Most of these classes have been already considered in \cite{HirataHespanha}. For the class of persistently exciting signals, see for instance \cite{PE-MCSS,Chitour-S,Sr-Ak} and references therein. Notice that all the classes in the above list are shift-invariant.

Rather than addressing the issues at stake for each class of switching signal given above, we develop a unifying framework which can also be applied to other classes. For that purpose, we adopt an axiomatic approach which singles out and exploits some useful common properties satisfied by the classes above.

\section{Adapted norms for switched linear systems with concatenable 
subfamilies}
\label{s:barabba}

We  consider in this section a switched linear system
\begin{equation}\label{SSA}
\dot x(t)=A(t)x(t)
\end{equation}
where $A$ belongs to a class $\mathcal{S}$ of measurable switching laws taking values in a bounded nonempty
set of  matrices $\mathscr{M}\subset M_n(\mathbb{R})$. 

A useful assumption on the family $\mathcal{S}$ that we are going to use in the following (which is satisfied by all the classes introduced in the previous section)  concerns its invariance by time-shift. 
\begin{itemize} %\begin{description} 
\item[{\bf A0}]{\bf (shift-invariance)} For every $A(\cdot)\in \mathcal{S}$ and every $t\geq 0$, the signal $A(t+\cdot)$ is in $\mathcal{S}$.
\end{itemize} %\end{description} 
Under Assumption {\bf A0} a convenient measure of the asymptotic behavior of \eqref{SSA} is the \emph{generalized spectral radius} (see, e.g., \cite{wirth-gsr})
\begin{equation}\label{eq:bohl}
\rho(\mathcal{S})=\limsup_{t\to +\infty}\sup_{A\in \mathcal{S}}\left\|\Phi_A(t,0)\right\|^{1/t}.
\end{equation}
Notice that, since $\mathscr{M}$ is bounded then $\rho(\mathcal{S})$ is finite.

As mentioned in introduction, our approach aims at extending the Barabanov norm construction (cf.~\cite{wirth-gsr}) beyond the class of signals with arbitrary switching. The main difficulty to do so
lies in the fact that the set of all flows $\Phi_A(s+t,s)$, for $A\in\mathcal{S}$ and $s,t\geq 0$, 
does not form a semigroup, since in general signals in $\mathcal{S}$ cannot be concatenated arbitrarily within $\mathcal{S}$.
A key object in what follows is then the identification of a subclass $\hS$ of $\mathcal{S}$, constructed by concatenating 
in an arbitrary way some signals defined on finite intervals. 
We then attach to $\hS$ a semigroup of fundamental matrices that captures the asymptotic behaviour of $\mathcal{S}$ if $\hS$ is large enough. 
 
\subsection{Concatenable subfamilies}\label{f:conca}

Consider a set $\mathcal{F}=\cup_{t\geq 0}\mathcal{F}_t$ with $\mathcal{F}_t\subset L_\infty([0,t],\mathscr{M})$, $t\in[0,\infty)$. Define 
\begin{align}\label{Shat}
\hS=\bigg\{&
A_1\ast A_2\ast\dots\ast A_k\ast\dots \mid \nonumber \\ 
&A_k \in \mathcal{F}_{t_k}\mbox{ for }k\in\N,\;\sum_{k\in\N}t_k=\infty
\bigg\}
\end{align}
and 
$\Mexp(\mathcal{F})=\cup_{t\ge0} \Mexp(\mathcal{F}_t)$, where, for every $t\ge0$,
$$\Mexp(\mathcal{F}_t)=\left\{ \Phi_A(t,0)\mid A\in \mathcal{F}_t\right\}.$$
Let, moreover,
 $$\mu(\mathcal{F})=\limsup_{t\to +\infty}\left(\sup\left\{\left\|R_t\right\|^{1/t}\mid R_t\in \Mexp(\mathcal{F}_{t})\right\}\right),
 $$
 with the convention that the quantity inside the parenthesis is equal to $0$ if $\mathcal{F}_{t}$ is empty.
Notice that $\mu(\mathcal{F})\leq  \rho(\hS)$,  but the converse is in general not guaranteed since the computation of $\rho(\hS)$ takes into account all intermediate instants between two concatenation times, unlike the one of $\mu(\mathcal{F})$.

We list below some useful  
 assumptions on the pair $(\mathcal{S},\mathcal{F})$ that will be exploited in the sequel.

%\begin{description}
\begin{itemize}
\item[{\bf A1}]{\bf (concatenability)} $\mathcal{F}_s\ast \mathcal{F}_t\subset \mathcal{F}_{s+t}$ for every $s,t\ge0$;
\item[{\bf A2}]  {\bf (irreducibility)} $\Mexp(\mathcal{F})$  is irreducible; 
\item[{\bf A3}] 
{\bf (fatness)} {$\hS\subset\mathcal{S}$ and} there exist two constants $C,\Delta\geq 0$ and a compact subset $\mathscr{K}$ of $\mathrm{GL}(n)$ 
such that for every $t\geq 0$  
and $A\in \mathcal{S}$,
 there exist $K\in \mathscr{K}$, $\hat t\in [t,t+\Delta]$, and  $R\in \Mexp(\mathcal{F}_{\hat t})$  such that
\begin{equation}\label{fat}
\left\| \Phi_A(t,0)\,KR^{-1}\right\|\leq C. 
\end{equation}
Moreover, if $A\in \hS$, one can take $K=\Id_n$ in \eqref{fat}. 
\end{itemize}
%\end{description}

\begin{remark}\label{conca0}
As a consequence of the definition of $\hS$, 
if $A\in \mathcal{F}$ and $B\in \hS$, then $A\ast B\in \hS$.
Moreover, by Assumption {\bf A1}, one has that 
$\Mexp(\mathcal{F}_s) \Mexp(\mathcal{F}_t)\subset \Mexp(\mathcal{F}_{s+t})$ for
every $s,t\geq 0$.
Hence $\hR$ is a semigroup and Assumption {\bf A2} above 
is equivalent to 
%\begin{equation}
%\label{A2bis}  
\[
\forall x\in \R^n\setminus\{0\},\mbox{ the linear 
span of $\hR x$ is equal to $\R^n$.}
\]
%\end{equation}
As in \cite{wirth-gsr}, one then says that $\hR$ is an irreducible semigroup. 
Note that \cite{wirth-gsr} considers special classes of irreducible semigroups which verify the additional assumption
\begin{center}
{\bf (decomposability)}~~$ \mathcal{F}_{s+t}=\mathcal{F}_s\ast \mathcal{F}_t$ for every $s,t\ge0.$
\end{center}
%%\begin{description}
%\begin{itemize}
%\item[\bf (decomposability)] $ \mathcal{F}_{s+t}=\mathcal{F}_s\ast \mathcal{F}_t$ for every $s,t\ge0.$
%\end{itemize}
%%\end{description}
The decomposability assumption trivially implies that $\Mexp(\mathcal{F}_s) \Mexp(\mathcal{F}_t)= \Mexp(\mathcal{F}_{s+t})$ for
every $s,t\geq 0$.
\end{remark}
 
\begin{remark}\label{clos}
Recall that the map $L_\infty([0,t],\mathscr{M})\ni A\mapsto \Phi_A(t,0)\in M_n(\R)$ is continuous with respect to the weak-$\star$ topology in $L_\infty([0,t],\mathscr{M})$ (see, for instance, \cite[Proposition 21]{PE-MCSS}). In particular if {\bf A3} holds true for 
$\mathcal{S}$ then it also holds true for 
the weak-$\star$ closure of $\mathcal{S}$.
\end{remark}

\begin{ml}\label{mu-lamb}
If {\bf A3} holds true then $\rho(\mathcal{S})=\mu(\mathcal{F})$. If moreover {\bf A0} and  {\bf A1} hold true then the family $\mathcal{F}$ verifies the following version of Fenichel's uniformity lemma: assume that there exists $m>0$ such that, for every sequence  $R_{j}\in \Mexp(\mathcal{F}_{t_j})$ with $t_j$ tending to infinity, one has $\lim_{j\to\infty} m^{-t_j} R_{j} = 0$. Then $\mu(\mathcal{F})<m$.
\end{ml}
%\IEEEproof 
\begin{IEEEproof}
The inequality $\rho(\mathcal{S})\geq\mu(\mathcal{F})$ is immediate. The opposite one readily comes from Assumption {\bf A3} and the definitions of $\rho(\mathcal{S})$ and $\mu(\mathcal{F})$. As regard the second part of the lemma, we can assume that $m=1$ by replacing if necessary $\mathscr {M}$ by the set $\mathscr {M} - \log (m) \mathrm{Id}_n$. Let $\mathcal{S}^\star$ be the  closure of $\mathcal{S}$ with respect to the weak-$\star$ topology induced by $L_\infty([0,\infty),\mathscr{M})$.
We first show that %,  %under the hypotheses on $\mathcal{F}$, one has that 
every trajectory associated with a switching signal in $\mathcal{S}^\star$ tends to zero. Assume  by contradiction that there exist $A\in\mathcal{S}^\star$, $x\in\R^n$, $\varepsilon>0$ and a sequence $t_j$ tending to infinity such that 
$$
\|\Phi_A(t_j,0)x\|\geq \varepsilon.
$$
By Remark~\ref{clos},  Assumption {\bf A3} actually extends to any switching signal in $\mathcal{S}^\star$. Hence, for every $j\geq 0$ applying Assumption {\bf A3} to the switching signal $A$ and the time $t_j$ yields the inequality 
$$\|\hat R_{j}K_j^{-1} x\|\geq \varepsilon/C$$
for some $\hat R_{j}\in  \Mexp(\mathcal{F}_{\hat{t}_j})$, with $\hat{t}_j\in[t_j,t_j+\Delta]$, and $K_j$ belonging to a given compact of $\mathrm{GL}(n)$. According to the hypotheses of the lemma, the left-hand side of the above inequality tends to $0$ as $j$ goes to infinity, which is  a contradiction.

Since the switching laws of $\mathcal{S}$ take values in the bounded set $\mathscr{M}$, one has that the class $\mathcal{S}^\star$ is weak-$\star$ compact. Recalling that $\mathcal{S}$ is shift-invariant by {\bf A0}, we notice that  all the assumptions of Fenichel's uniformity lemma are verified by the standard linear flow defined on $\R^n\times \mathcal{S}^\star$ (cf.~\cite{colonius-dynC}). Therefore the convergence of trajectories of $\mathcal{S}$ to $0$ is uniformly exponential, i.e., $\rho(\mathcal{S})<1$. By the first part of the lemma we deduce that $\mu(\mathcal{F})<1$. 
\end{IEEEproof}

We associate with each class of switching signals considered in the previous section a corresponding concatenable subfamily as listed below: 
\begin{itemize}
\item $\mathcal{F}^{\rm arb}(\mathscr{M})$: arbitrarily switching signals on finite intervals;
\item  $\mathcal{F}^{\rm pc}(\mathscr{M})$: piecewise constant signals on finite intervals;
\item $\mathcal{F}^{{\rm d},\tau}(\mathscr{M})$: piecewise constant signals on finite intervals with dwell-time $\tau$ {and such that the first and last subintervals on which the signal is constant have length at least  $\tau$} (notice that $\mathcal{F}_t^{{\rm d},\tau}(\mathscr{M})=\emptyset$ for $t<\tau$); 
\item $\mathcal{F}^{{\rm av-d},\tau,N_0}(\mathscr{M})$: piecewise constant signals on finite intervals  satisfying the 
$(\tau,N_0)$ average dwell-time condition and such that the {first and} last subintervals on which the signal is constant {have} length at least  $N_0\tau$;
 \item $\mathcal{F}^{{\rm pe},T,\mu}(\mathscr{M})$: $(T,\mu)$-persistently exciting signals on finite intervals $[0,t]$ for which $t\geq T-\mu$ and the signal is constantly equal to $M_1$ on $[t-T+\mu,t]$, where  we recall that $\mathscr{M}=\{M_{{\delta}}=(1-\delta)M_0+{\delta} M_1\mid {\delta}\in [0,1]\}$.
\end{itemize}
For the classes $\mathcal{S}^{{\rm lip},L}(\mathscr{M})$ and $\mathcal{S}^{{\rm BV},T,\mu}(\mathscr{M})$ 
we fix some $\widebar M \in \mathscr{M}$ and we define
\begin{itemize}
 \item $\mathcal{F}^{{\rm lip},L}(\mathscr{M})$: 
Lipschitz signals on finite interval starting and ending at $\widebar M$;
 \item $\mathcal{F}^{{\rm BV},T,\mu}(\mathscr{M})$: $(T,\nu)$-BV signals on finite intervals $[0,t]$, $t\geq T$,  starting and ending at $\widebar M$ and constant on $[t-T,t]$. 
\end{itemize}
With these choices of $\mathcal{F}$ Assumption ${\bf A1}$ is automatically satisfied. %(with {\bf B1} holding only for $\mathcal{F}^{\rm arb}(\mathscr{M})$ and $\mathcal{F}^{\rm pc}(\mathscr{M})$). 
To address the validity of Assumption ${\bf A2}$ for the previous classes of signals, we further introduce the following assumption on the set $\mathcal{F}$, which essentially says that the flow corresponding to any element in 
$\mathcal{F}^{\rm pc}(\mathscr{M})$ can be approached in a suitable sense by an analytic deformation of flows corresponding to elements in $\mathcal{F}$. 
\begin{itemize} %\begin{description} 
\item[{\bf A4}]{\bf (Analytic propagation)}  For every $t,\varepsilon>0$ and $A\in \mathcal{F}_t^{\rm pc}(\mathscr{M})$, there exists a path $(0,1]\ni \lambda\mapsto (t_\lambda,A_{\lambda})\in (0,\infty)\times \mathcal{F}_{t_\lambda}^{\rm arb}(\mathscr{M})$  such that 
\begin{itemize}
\item $\lambda\mapsto \Phi_{A_{\lambda}}(t_\lambda,0)$ is analytic; 
\item $\Vert  \Phi_{A_1}(t_1,0)- \Phi_A(t,0)\Vert\leq \varepsilon$;
\item the set $\{\lambda\in (0,1]\mid A_{\lambda}\in  \mathcal{F}\}$ has nonempty interior.
\end{itemize}
\end{itemize} %\end{description}
The relation between Assumptions {\bf A2} and {\bf A4} is clarified in the following proposition.
\begin{proposition}\label{prop:A2}
%Consider a set $\mathcal{F}=\cup_{t\geq 0}\mathcal{F}_t$ with $\mathcal{F}_t\subset L_\infty([0,t],\mathscr{M})$, $t\ge0$, where 
Let $\mathscr{M}$ be irreducible. Then Assumption {\bf A4} implies Assumption {\bf A2}.
\end{proposition}
%\IEEEproof 
\begin{IEEEproof}
Let $x,z\in \R^n\setminus\{0\}$. We have to prove that there exists $R\in \hR$ such that 
$z^T Rx\ne 0$. 
Since $\mathscr{M}$ is irreducible, it follows from \cite[Lemma 3.1]{wirth-gsr} that there exist $t>0$ and $A$ in $\mathcal{F}_t^{\rm pc}(\mathscr{M})$ such that 
%\begin{equation}\label{nonz}
\[
z^T \Phi_A(t,0) x\ne 0.
\]
%\end{equation}
Consider the path  
$\lambda\mapsto A_{\lambda}$ provided by Assumption {\bf A4}.
The function  
$$
\lambda \mapsto z^T \Phi_{ A_\lambda}(t_\lambda,0) x
$$ 
is analytic and not identically equal to zero. It therefore vanishes at isolated values of $\lambda$, whence the conclusion with $R$ of the type $\Phi_{ A_\lambda}(t_\lambda,0)$. 
\end{IEEEproof}

In the following proposition we establish the
validity of Assumptions {\bf A0}, {\bf A1}, {\bf A3} and {\bf A4} for the classes 
of switching signals introduced in Section \ref{s:classes} and their corresponding concatenable subfamilies.

\begin{proposition}\label{lem:A4}
Let $\mathcal{S}$ be one of the classes introduced in Section~\ref{s:classes} with corresponding 
$\mathcal{F}$ as above. 
If $\mathcal{S}=\mathcal{S}^{\mathrm{lip},L}(\mathscr{M})$ or $\mathcal{S}=\mathcal{S}^{{\rm BV},T,\nu}(\mathscr{M})$, assume moreover that $\mathscr{M}$ is star-shaped, that is, 
there exists $\hat M\in \mathscr{M}$ such that for any other $M\in\mathscr{M}$ the segment between $\hat M$ and $M$ is contained in $\mathscr{M}$. 
Then Assumptions {\bf A0}, {\bf A1}, {\bf A3} and {\bf A4} hold true. 
\end{proposition}
%\IEEEproof 
\begin{IEEEproof}
As already noticed, Assumptions {\bf A0} and {\bf A1} %and {\bf A3} 
are satisfied.

Concerning Assumption ${\bf A3}$, notice that every restriction $A|_{[0,t]}$ of a signal in one of the classes $\mathcal{S}$ introduced in Section~\ref{s:classes} can be extended to a signal ${A_1}\ast A|_{[0,t]}\ast A_2$ in the corresponding class {$\mathcal{F}$}, with $A_j:[0,t_j]\to\mathscr{M}$, $j=1,2$, and $t_1,t_2\leq t_*$ for some $t_*$ uniform with respect to 
$A\in\mathcal{S}$ and $t\ge0$. 
Moreover, if $A\in\mathcal{F}_\infty$ then $t_1$ can be taken equal to zero.

Let us now prove Assumption ${\bf A4}$. 
Take $t,\varepsilon>0$ and $A\in \mathcal{F}_t^{\rm pc}(\mathscr{M})$. 
For the cases $\mathcal{S}=\mathcal{S}^{\mathrm{d},\tau}(\mathscr{M})$ and $\mathcal{S}=\mathcal{S}^{\mathrm{av}-\mathrm{d},\tau,N_0}(\mathscr{M})$ one can find the path $\lambda\mapsto A_{\lambda}$ 
as follows. For every $\delta>0$, consider the time-reparameterized signal 
$A(\delta\, \cdot)\in \mathcal{F}_{t/\delta}^{\rm pc}(\mathscr{M})$. Then $A(\delta\, \cdot)\in \mathcal{F}$ for $\delta$ small enough and the function  
$\delta \mapsto \Phi_{A(\delta\,\cdot)}(t/\delta,0)$ is analytic. Indeed, the function $(0,\infty)\ni\delta\mapsto \Phi_{A(\delta\,\cdot)}(t_0/\delta,0)=\Phi_{ A(\cdot)/\delta}(t_0,0)$ is analytic, since the Volterra series associated with this flow define an analytic function of $\delta$.

We conclude by taking $A_{\lambda}(\cdot)=A(\lambda\,\cdot)$. 
For the cases $\mathcal{S}=\mathcal{S}^{\mathrm{pe},T,\mu}(\mathscr{M})$ and $\mathcal{S}=\mathcal{S}^{\mathrm{BV},T,\mu}(\mathscr{M})$, the previous construction can be modified as follows. First notice that, up to adding some short intervals on which $A$ is constant we can assume that
\begin{itemize} %\begin{description} 
%\item[Case $\mathcal{S}=\mathcal{S}^{\mathrm{pe},T,\mu}(\mathscr{M})$:] 
\item Case $\mathcal{S}=\mathcal{S}^{\mathrm{pe},T,\mu}(\mathscr{M})$: the value of $A$ on the last  interval on which it is constant is $M_1$;
%\item[Case $\mathcal{S}=\mathcal{S}^{\mathrm{BV},T,\nu}(\mathscr{M})$:]
\item Case $\mathcal{S}=\mathcal{S}^{\mathrm{BV},T,\nu}(\mathscr{M})$: the distance between two subsequent values of $A$ is smaller than $\nu$, and the first and last values of $A$ are both equal to $\widebar  M$.
\end{itemize} %\end{description} 
The modification of $A$ can be taken so that the corresponding variation of $\Phi_A(t,0)$  is smaller than $\varepsilon$. The argument now works as before in the case $\mathcal{S}=\mathcal{S}^{\mathrm{BV},T,\nu}(\mathscr{M})$. In the case $\mathcal{S}=\mathcal{S}^{\mathrm{pe},T,\mu}(\mathscr{M})$ one can take $A_{\lambda}(s) =A( s/ \lambda)$ on $[0,\lambda t_0] $ and $A_\lambda(s)=M_1$ on $[\lambda t_0,t]$, where $t_0$ is such that $A|_{[t_0,t]}\equiv M_1$.

We are left to discuss the case $\mathcal{S}=\mathcal{S}^{\mathrm{lip},L}(\mathscr{M})$.
Similarly to what is done above, we can first assume that $A$ is equal to $\widebar M$ on the first and last interval on which it is constant. Then we modify $A$ by adding, 
 at each switching time, a Lipschitz continuous arc defined on a small time interval and bridging the discontinuity (one may take a Lipschitz continuous arc whose graph is the union of two segments joining at $\hat M$).
These modifications can be done  while keeping   the corresponding variation of $\Phi_A(t,0)$ smaller than $\varepsilon$. 
By 
a time-reparameterization $\lambda\mapsto A(\lambda\,\cdot)$ we can lower the Lipschitz constant and complete the proof as above. 
\end{IEEEproof}

\begin{remark}
In Proposition~\ref{lem:A4}, the hypothesis on the star-shapedness  of $\mathscr{M}$ can be replaced by some weaker one. For instance we could assume:
\begin{itemize}
\item if $\mathcal{S}=\mathcal{S}^{\mathrm{lip},L}(\mathscr{M})$ then 
there exists $C>0$ such that  
every two distinct points of $\mathscr{M}$ can be connected by a Lipschitz-continuous curve lying in 
$\mathscr{M}$ of length smaller than $C$;
\item if $\mathcal{S}=\mathcal{S}^{{\rm BV},T,\nu}(\mathscr{M})$ then for every $M_0,M_1\in  \mathscr{M}$, there exists a finite sequence of points in $\mathscr{M}$ whose first and last elements are $M_0$ and $M_1$, respectively, and such that the distance between two subsequent elements is smaller than $\nu$.
\end{itemize}
\end{remark}

\subsection{Quasi-Barabanov semigroups}
 The main goal  of the section is to prove the result below.
\begin{Theorem}\label{anticipazioni}
Let  $(\mathcal{S},\mathcal{F})$ satisfy Assumptions~{\bf A0--A3}.
Then there exists 
a constant $C \ge1$ such that for any $x_0\in \R^n\setminus\{0\}$ there exists
a trajectory  $x:t\mapsto \Phi_A(t,0) x_0$ with  $A$ belonging to the weak-$\star$ closure of 
$\hS$
such that, for every $t\ge0$, 
\[
\frac1{C} \rho(\mathcal{S})^t\|x_0\|\le \|x(t)\|\le C \rho(\mathcal{S})^t\|x_0\|.
\]
\end{Theorem}
For that purpose, we first need the following definitions.

\begin{Definition} 
Let $\mathscr{M}$,  $\mathcal{F}$,  $\hS$, and $\hR$ be as in the previous section.
%\begin{enumerate}
%\item 

We say that $\hR$ is a  \emph{quasi-extremal semigroup} if there exists $C_{\mathrm{qe}}>0$ such that, for every $t\geq 0$ and 
$R\in  \Mexp(\mathcal{F}_t)$, one has
\begin{equation}\label{e:qest}
\|R\|\leq C_{\mathrm{qe}} \mu(\mathcal{F})^t.
\end{equation}
Moreover, a quasi-extremal semigroup $\hR$ is said to be \emph{extremal} if there exists a norm 
$w$ on $\R^n$ such that the induced matrix norm $\|\cdot\|_w$ satisfies, for every $t\ge 0$ and $R\in \Mexp(\mathcal{F}_t)$, 
\begin{equation}\label{e:qest-1}
\|R\|_w\le \mu(\mathcal{F})^t.
\end{equation}
A norm $w$ satisfying Eq.~\eqref{e:qest-1} is said to be  \emph{extremal} for $\hR$.
A quasi-extremal semigroup  $\hR$ is said to be  \emph{quasi-Barabanov}  if there exists $C_{\mathrm{qb}}>0$ such that 
for  every $x\in \R^n$ and $t\geq 0$ there exist $t'\geq t$ and  $R\in \Mexp(\mathcal{F}_{t'})$ 
such that
\begin{equation}\label{e:qb}
\|Rx\|\geq C_{\mathrm{qb}} \mu(\mathcal{F})^{t'}\|x\|.
\end{equation}
%\item
Let  $\hR$ be a quasi-extremal semigroup. A  trajectory $x:t\mapsto \Phi_A(t,0) x_0$ with $x_0\ne 0$ and 
$A$ belonging to the weak-$\star$ closure of 
$\hS$ 
is said to be \emph{
quasi-extremal with constant   $C_{\mathrm{qx}} \ge1$}  if for every $t\ge0$, 
%\begin{equation}\label{Cext}
\[  
  \frac1{C_{\mathrm{qx}}} \mu(\mathcal{F})^t\|x_0\|\le \|x(t)\|\le C_{\mathrm{qx}} \mu(\mathcal{F})^t\|x_0\|.
\]
%  \end{equation}
%\end{enumerate}
\end{Definition}
 
The notion of generalized spectral radius for a quasi-Barabanov semigroup $\hR$ is actually equivalent to the following adaptation of the definition of  maximal Lyapunov exponent
 \begin{align*}
 \lambda(\mathcal{F})=\sup\bigg\{& \limsup_{k\to \infty}\frac{\log\|R_{t_k}\cdots R_{t_1}\|}{t_1+\dots+t_k}\mid \\
& R_{t_k}\in\Mexp(\mathcal{F}_{t_k}) \mbox{ for $k\in\N$},\;\sum_{k\in\N}t_k=\infty\ \bigg\},
 \end{align*}
as stated below.
\begin{proposition}\label{p-Lyapunov}
Let $\mathcal{F}$ be a family of switching laws satisfying the concatenability condition {\bf A1} and assume that $\hR$ is a quasi-Barabanov semigroup. Then $ \mu(\mathcal{F})=e^{\lambda(\mathcal{F})}$.
\end{proposition}
%\IEEEproof %
\begin{IEEEproof}
The inequality $ \mu(\mathcal{F})\geq e^{\lambda(\mathcal{F})}$ easily comes from \eqref{e:qest}. In order to show the opposite inequality one observes that, by \eqref{e:qb}, for any $x_0\in\R^n$ and $i\in\N$, there exist $t_i\geq i$ and $R_{t_i}\in \Mexp(\mathcal{F}_{t_i})$ such that  
\[
\|R_{t_k}\cdots R_{t_1}x_0\| \geq (C_{\mathrm{qb}})^k \mu(\mathcal{F})^{t_1+\dots+t_k} \|x_0\|.
\]
In particular
\begin{align*}
\limsup_{k\to\infty}\frac{\log\|R_{t_k}\cdots R_{t_1}\|}{t_1+\dots+t_k} &\geq -\!\lim_{k\to\infty} \!\frac{2|\log(C_{\mathrm{qb}})|}{k+1} + \log(\mu(\mathcal{F}))\\
&=\log(\mu(\mathcal{F})),
\end{align*}
which concludes the proof.
\end{IEEEproof}
\begin{ml}\label{qB-qE}
Let Assumptions {\bf A0}, {\bf A1} 
 and {\bf A3} be satisfied. 
If $\hR$ is a quasi-Barabanov semigroup then 
there exists $C_{\mathrm{qx}}\geq 1$ such that 
any nonzero point of $\R^n$
is the initial condition of a quasi-extremal trajectory with constant $C_{\mathrm{qx}}$.
\end{ml}
%\IEEEproof %
\begin{IEEEproof}
Let $C_{\rm qe}$, $C_{\rm qb}$ be as in \eqref{e:qest}, \eqref{e:qb} and $C$, $\mathscr{K}$  as in {\bf A3}.
Let $\kappa\ge 1$ verify $\|K\|,\|K^{-1}\|
\leq 
\kappa$ for every 
$ K\in\mathscr{K}$.

Fix $x_0\ne 0$.
There exists an increasing sequence of times  $t_k$ going to infinity and signals $A_k$ in 
$\mathcal{F}_{t_k}$ such that 
$\|R_k x_0\|\geq C_{\mathrm{qb}} \mu(\mathcal{F})^{t_k}\|x_0\|$, where $R_k=\Phi_{A_k}(t_k,0)$.

We now claim that,
for every $k\in \N$ and every $s\in [0,t_k]$, one has
\begin{equation}\label{zor}
 \left\|\Phi_{A_k}(s,0) x_0\right\|\geq \frac1{C_0} \, \mu(\mathcal{F})^s\|x_0\|
 \end{equation}
 for some constant $C_0>0$ independent of $x_0$ and $s$. 
Indeed, 
because of Assumption {\bf A0} and 
applying \eqref{fat} 
one gets that
$\Phi_{A_k}(t_k,s)=M R K^{-1}$ with $\|M\|\le C$, $K\in\mathscr{K}$, and $R$ in 
$\Mexp(\mathcal{F}_{t_k-s+\delta})$, for some $\delta\in[0,\Delta]$. 
One can therefore write $R_k=M R K^{-1}\Phi_{A_k}(s,0)$. It follows that
\begin{align*}
 \|R_k x_0\|&\leq  \|M\| \left\|R K^{-1} \Phi_{A_k}(s,0) x_0\right\|\\
& \le C C_{\mathrm{qe}} \mu(\mathcal{F})^{t_k-s+\delta}\left\|K ^{-1} \Phi_{A_k}(s,0) x_0\right\|\\
 &\le \kappa CC_{\mathrm{qe}} \mu(\mathcal{F})^{t_k-s+\delta}\left\|\Phi_{A_k}(s,0) x_0\right\|.
 \end{align*}
 On the other hand 
 $\|R_k x_0\|\ge C_{\mathrm{qb}} \mu(\mathcal{F})^{t_k} \|x_0\|$, which 
 proves \eqref{zor} with
 $$C_0=\frac{\kappa C C_{\rm qe}\max(1,\mu(\mathcal{F})^{\Delta})
 }{C_{\mathrm{qb}}} .$$ 

Notice that each $A_k$ is the restriction on $[0,t_k]$ of a signal $B_k$ in 
$\hS$. 
Up to a subsequence we can assume that $B_k$ weak-$\star$ converges to some $B_\star$ 
in the weak-$\star$ closure of  
$\hS$.
Passing to the limit in \eqref{zor} 
we deduce that for every $s\geq 0$, 
\begin{equation}\label{left0} 
\left\|\Phi_{B_\star}(s,0) x_0\right\|\geq \frac 1{C_0} \, \mu(\mathcal{F})^s\|x_0\|.
\end{equation}

We next prove that 
there exists $C_1>0$ such that
for every $s\geq 0$, $x_0\in\R^n$ and   $B$ in the weak-$\star$ closure of $\hS$, it holds
\begin{equation}\label{right0} 
\left\|\Phi_{B}(s,0) x_0\right\|\leq C_1 \, \mu(\mathcal{F})^s\|x_0\|.
\end{equation}
For that purpose, 
consider a sequence $B_k$ in $\hS$ weak-$\star$ converging to $B$. 
Applying \eqref{fat}, Remark~\ref{clos} 
and a compactness argument, we get that $\Phi_{B}(s,0)=MR K^{-1}$ with $\|M\|\le C$, $K\in\mathscr{K}$ and $R$ in the closure of $\cup_{\delta\in[0,\Delta]}\Mexp(\mathcal{F}_{s+\delta})$. 
 It follows that
$$
\left\|\Phi_{B}(s,0) x_0\right\|\leq  \|M\|  \|RK^{-1} x_0\|
\le C_1\mu(\mathcal{F})^{s}\|x_0\|,
 $$
where $C_1=\kappa C C_{\mathrm{qe}}  \max(1,\mu(\mathcal{F})^{\Delta})$
, proving \eqref{right0}. 
Together with   \eqref{left0}, this concludes the proof of the lemma with $C_{\mathrm{qx}}=\max(C_0,C_1)$.
\end{IEEEproof}

Set 
\begin{align*}
\Rinfty=\{R\mid\ & \exists t_k\to\infty,\ R_{k}\in \Mexp(\mathcal{F}_{t_k})\\
&\mbox{such that }
\mu(\mathcal{F})^{-t_k}R_{k}\to R\}.
\end{align*}
The following result is the counterpart of~\cite[Proposition 3.2]{wirth-gsr} in our setting.

\begin{mpr}\label{prop1}
Let  $(\mathcal{S},\mathcal{F})$ satisfy Assumptions~{\bf A0--A3} and define 
$\Rinfty$ as above. 
Then 
\begin{itemize} %\begin{description} 
\item[$(i)$] $\Rinfty$ is compact and nonempty, $\Rinfty\neq \{0\}$;
\item[$(ii)$] $\Rinfty$ is a semigroup;
\item[$(iii)$] for every $t\geq 0$, $T\in \Mexp(\mathcal{F}_t)$ and $S\in \Rinfty$, both $\mu(\mathcal{F})^{-t}TS$ and $\mu(\mathcal{F})^{-t}ST$ belong to
$\Rinfty$;
\item[$(iv)$] $\Rinfty$ is irreducible.
\end{itemize} %\end{description} 

\end{mpr}
\begin{IEEEproof}
To prove the proposition one follows exactly the arguments provided in the proof of \cite[Proposition 3.2]{wirth-gsr} except for the fact that $\Rinfty\neq \{0\}$. In our setting this result is easily proved by using Lemma~\ref{mu-lamb} and the definition of $\Rinfty$.
\end{IEEEproof}

The following result can be proven as in \cite[Lemma 3.4]{wirth-gsr}.
\begin{mpr}\label{prop2} 
Let $(\mathcal{S},\mathcal{F})$ satisfy Assumptions~{\bf A0--A3} and   define 
$\Rinfty$ as above. 
Let $\hat v:\R^n\rightarrow (0,\infty)$ be 
defined as
%\beq\label{nBar0}
\[
\hat v(x)=\max_{R\in \Rinfty}\Vert Rx\Vert.
\]
%\eeq
Then $\hat v$ is an extremal norm for $\hR$.
\end{mpr}

\begin{remark}
If Assumption {\bf A1} is replaced by the stronger decomposability assumption (see Remark~\ref{conca0}), then $\hat v$ is a Barabanov norm (see e.g. \cite{wirth-gsr}).  
\end{remark}

We have the following result.
\begin{mpr}\label{prop3} 
Let $(\mathcal{S},\mathcal{F})$ satisfy Assumptions~{\bf A0--A3}.
Then 
$\hR$ is an extremal and quasi-Barabanov semigroup. 
\end{mpr}
\begin{IEEEproof} 
 The fact that $\hR$ is an extremal semigroup readily comes from Proposition~\ref{prop2}. In order to show the second part of the statement, let us
consider $\kappa\geq 1$ 
as in the proof of Lemma~\ref{qB-qE}.
Without loss of generality, let us also assume that $\kappa$ satisfies
$$
\frac1\kappa \Vert x\Vert \leq \hat{v}(x)\leq \kappa \Vert x\Vert\hbox{ for all }x\in\mathbb{R}^n
$$
and that $\|M\|\leq C$ implies that $\|M\|_{\hat v}\le\kappa$, where $C$ is as in {\bf A3}.

For every $x_0\in\mathbb{R}^n$, there exists a sequence $(R_k)_{k\in \N}$ so that $R_k\in \Mexp(\mathcal{F}_{t_k})$  with $\lim_{k\rightarrow\infty}t_k=\infty$ and    $\hat v(x_0)=\lim_{k\rightarrow\infty}\mu(\mathcal{F})^{-t_k}\Vert R_kx_0\Vert$. For every $k\in\N$, let 
$A_k\in\mathcal{F}_{t_k}$ 
be such that $R_k=\Phi_{A_k}(t_k,0)$.

Fix now $s\geq 0$. For every $k\in \N$ such that $t_k\geq s$, 
we deduce from  Assumption~{\bf A3} that $\Phi_{A_k}(s,0)=M_1^{(k)}Q_1^{(k)}$ and $\Phi_{A_k}(t_k,s)=M_2^{(k)}Q_2^{(k)}(K_2^{(k)})^{-1}$ with $\|M_i^{(k)}\|\le C$, $K_2^{(k)}\in\mathscr{K}$, $Q_1^{(k)}\in  \Mexp(\mathcal{F}_{s+\delta_1^{(k)}})$ and  $Q_2^{(k)}\in \Mexp(\mathcal{F}_{t_k-s+\delta_2^{(k)}})$, where  $\delta_i^{(k)}\in[0,\Delta]$, for $i=1,2$.
One can therefore write $R_k=M_2^{(k)}Q_2^{(k)}(K_2^{(k)})^{-1}M_1^{(k)}Q_1^{(k)}$. Using the extremality of $\hat{v}$ and the requirements imposed on $\kappa$, it follows that
\begin{align*}
 \mu(\mathcal{F})^{-t_k} \Vert R_kx_0\Vert&\leq \kappa \mu(\mathcal{F})^{-t_k}\hat v(R_k x_0)\\
& \leq \kappa^6 \mu(\mathcal{F})^{-s+\delta_2^{(k)}}
 \hat v(Q_1^{(k)}x_0).
\end{align*}
Taking limits as $k$ tends to infinity and up to subsequences, 
one gets 
$$
\hat{v}(x_0)\leq \kappa^6 \mu(\mathcal{F})^{-s+\delta_2}\hat{v}(Q_1 x_0),
$$
where $Q_1$ is in the closure of $\cup_{\delta\in[0,\Delta]}\Mexp(\mathcal{F}_{s+\delta})$, and
$\delta_2$ belongs to $[0,\Delta]$.
In particular, there exist $\hat Q\in \Mexp(\mathcal{F}_{s+\hat \delta})$ for some $\hat \delta\in[0,\Delta]$ such that
$$
\hat{v}(x_0)\leq 2\kappa^6 \mu(\mathcal{F})^{-s+\delta_2}\hat{v}(\hat Q x_0).
$$
One deduces that $\hat{v}(\hat Q x_0)\geq C_0 \mu(\mathcal{F})^{s+\hat \delta}\hat{v}(x_0)$ with  $C_0=\frac{\min(1,\mu(\mathcal{F})^{-2\Delta})}{2\kappa^6}$.
 \end{IEEEproof}

As a consequence of Lemma~\ref{qB-qE}, we have the following corollary.
\begin{corollary}\label{qui-c}
Let $(\mathcal{S},\mathcal{F})$ satisfy Assumptions~{\bf A0--A3}. 
Then
there exists $C_{\mathrm{qx}}\geq 1$ such that
any nonzero point of $\R^n$ is the initial condition of a 
quasi-extremal trajectory with constant $C_{\rm qx}$.
\end{corollary}

Theorem~\ref{anticipazioni} follows directly from Corollary~\ref{qui-c} and Lemma~\ref{mu-lamb}.

Note that, as a consequence of  {\bf A3} and the above results, the right-hand side inequality in the statement of Theorem~\ref{anticipazioni} holds true for any trajectory associated with a signal in $\mathcal{S}$, up to adapting the constant $C$.
\begin{remark}
As a consequence of Proposition~\ref{p-Lyapunov}, Proposition~\ref{prop3} and Lemma~\ref{mu-lamb} one easily deduces that $\rho(\mathcal{S})=e^{\lambda(\mathcal{S})}$, where 
\[\lambda(\mathcal{S})=\sup_{A\in \mathcal{S}}\limsup_{t\to +\infty}\frac{\log\|\Phi_{A}(t,0)\|}{t}\]
is the \emph{maximal Lyapunov exponent} associated with the family $\mathcal{S}$. This result was already obtained in \cite{Wirth-SIAM} if the class $\mathcal{S}$ is assumed to be weak-$\star$ closed. 
\end{remark}

 \section{$L_2$-gain for switched linear control systems}
 \label{s:L2}

Consider a switched linear control system of the type 
 \begin{equation}
% \begin{split}
 \dot x(t)=A(t)x(t)+B(t)u(t),\ y(t)=C(t)x(t),\label{LSCS-D}
%\end{split} 
\end{equation}
 where $(x,u,y)\in \mathbb{R}^n\times\mathbb{R}^m\times\mathbb{R}^p$ and $(A,B,C)$ belongs to a class $\mathcal{T}$ of measurable switching laws taking values in a bounded set of triples of matrices $\mathscr{M}\subset M_n(\mathbb{R})\times  M_{n, m}(\mathbb{R})\times  M_{p,n}(\mathbb{R})$. 
We denote by $\pi_A$ and $\pi_{A\times B}$ the projections from $M_n(\mathbb{R})\times  M_{n, m}(\mathbb{R})\times  M_{p,n}(\mathbb{R})$ to its first and first two factors, respectively. 
  We set 
  \begin{align*}
  &\mathscr{M}_A=\pi_A(\mathscr{M})=\{A\mid \exists (A,B,C)\in \mathscr{M}\}, \\
   &\mathcal{T}_A=\pi_A(\mathcal{T})=\{A\mid \exists (A,B,C)\in \mathcal{T}\}\end{align*}
  and we define similarly $\mathscr{M}_B$, $\mathscr{M}_C$, $\mathcal{T}_B$ and $\mathcal{T}_C$.

 In the sequel, we also assume that $\mathcal{T}$ 
 contains a subset $\mathcal{G}_\infty$ made of concatenations of signals in a family $\mathcal{G}=\cup_{t\geq 0}\mathcal{G}_t$ as in Eq.~\eqref{Shat}.
Useful properties on  $(\mathcal{T},\mathcal{G})$ 
 are the following:
%\begin{description}
\begin{itemize}
\item[{\bf T1}] $(\mathcal{T}_A,\mathcal{G}_A)$ satisfies Assumptions {\bf A0} and {\bf A1}, where $\mathcal{G}_A=\pi_A(\mathcal{G})$. 
\item[{\bf T2}] $(\mathcal{T}_A,\mathcal{G}_A)$ satisfies Assumption {\bf A4}. Moreover, for every $t^*,\varepsilon>0$ and $(A,B,C)\in \mathcal{F}_t^{\rm pc}(\mathscr{M})$, there exists a path $[0,1]\ni \lambda\mapsto (t_\lambda,A_\lambda,B_\lambda,C_\lambda)\in (0,\infty)\times \mathcal{F}_{t_\lambda}^{\rm arb}(\mathscr{M})$ such that: 
\begin{itemize}
\item $\lambda\mapsto (W_\lambda^c(t_\lambda),W_\lambda^o(t_\lambda))$ is analytic, where $W_\lambda^c(t_\lambda)$ and $W_\lambda^o(t_\lambda)$ denote, respectively, the controllability and observability Gramians in time $t_\lambda$ associated with $\dot x(t) =A_\lambda(t)x(t)+B_\lambda(t)u(t)$, $y(t)=C_\lambda(t)x(t)$;
\item $\Vert W^c_1(t_1)- W^c(t^*)\Vert \leq \varepsilon$ and $\Vert W^o_1(t_1)- W^o(t^*)\Vert \leq \varepsilon$, where $W^c(t^*)$ and $W^o(t^*)$ denote, respectively, the controllability and observability Gramians in time $t^*$ associated with $\dot x(t) =A(t)x(t)+B(t)u(t)$, $y(t)=C(t)x(t)$;
\item the set $\{\lambda\in [0,1]\mid (A_\lambda,B_\lambda,C_\lambda)\in  \mathcal{G}\}$ has nonempty interior.
\end{itemize}
\end{itemize}
%\end{description}
A trivial adaptation of the proof of Proposition~\ref{lem:A4} yields the following result.
\begin{ml}\label{T1-T2} Let $\mathcal{S}$ be one of the classes introduced in Section~\ref{s:classes} with corresponding 
$\mathcal{F}$. 
Assume moreover that if $\mathcal{S}=\mathcal{S}^{\mathrm{lip},L}(\mathscr{M})$ 
or %then $\mathscr{M}$ is rectifiably pathwise connected and that if 
$\mathcal{S}=\mathcal{S}^{{\rm BV},T,\nu}(\mathscr{M})$ then 
%for every $M_0,M_1\in  \mathscr{M}$, there exists a finite sequence of points in $\mathscr{M}$ whose first and last elements are $M_0$ and $M_1$, respectively, and such that the distance between two subsequent elements is smaller than $\nu$.
$\mathscr{M}$ is star-shaped.
Then $(\mathcal{S},\mathcal{F})$ satisfies Assumptions {\bf T1} and {\bf T2}. 
\end{ml}

\subsection{Minimal realization for switched linear control systems}
We start by giving the following definitions. 
 
\begin{Definition}\label{rea-obs}
\begin{enumerate}
\item A point $x\in\mathbb{R}^n$ is $\mathcal{G}$-reachable  for the switched linear control system \eqref{LSCS-D}  if there exist $t\ge0$, a switching law $(A,B,C)\in\mathcal{G}_t$ and an input $u\in L_2$ such that the corresponding trajectory $x_u$ starting at $0$ reaches $x$ in time $t$. 
The \emph{reachable set} $\mathfrak{R}(\mathcal{G})$  is the set of 
all $\mathcal{G}$-reachable 
points.   System \eqref{LSCS-D} is said to be $\mathcal{G}$-controllable if $\mathfrak{R}(\mathcal{G})=\mathbb{R}^n$.
\item
 A point $x\in\mathbb{R}^n$ is $\mathcal{G}$-observable for the switched linear control system \eqref{LSCS-D} if there exist $t\ge0$ and a switching law $(A,B,C)\in\mathcal{G}_t$ such that the trajectory $x_0$ associated with the zero input and starting at $x$ gives rise to an output $y$ verifying $y(t)\neq 0$. The \emph{observable set} $\mathfrak{O}(\mathcal{G})$  is the set of all $\mathcal{G}$-observable points. System \eqref{LSCS-D} is said to be $\mathcal{G}$-observable if $\mathfrak{O}(\mathcal{G})=\mathbb{R}^n$.
\end{enumerate}
\end{Definition}

It is not immediate from their definitions that the reachable and observable sets $\mathfrak{R}(\mathcal{G})$  and $\mathfrak{O}(\mathcal{G})$ are linear subspaces. It has been shown in \cite{Sun-Ge-Lee} that this is the case if $\mathcal{T}=\mathcal{G}_\infty=\mathcal{S}^{\rm arb}(\mathscr{M})$, where 
$\mathscr{M}=\{(A_1,B_1,C_1),\dots,(A_k,B_k,C_k)\}$ with $k$ a positive integer. In addition, it is proved in the same reference that the state space admits a direct sum decomposition into a controllable and an uncontrollable part for the switched linear control system exactly as in the unswitched situation. 
More precisely, there exists a direct sum decomposition of the state space $\mathbb{R}^n=\mathfrak{R}(\mathcal{G})\oplus E$ and an invertible $n\times n$ matrix $P$ such that, if $r=\dim \mathfrak{R}(\mathcal{G})$ and $P^{-1}x=(x^c\ x^u)$ with $x^c\in\mathfrak{R}(\mathcal{G})$ and $x^u\in E$, 
one has for $1\leq i\leq k$,
$$
P^{-1}\!A_iP=\begin{pmatrix}A_i^c & \ast\\0 & A_i^u\end{pmatrix}\!,\ \
P^{-1}\!B_i=\begin{pmatrix}B_i^c\\0\end{pmatrix}\!,\ \
C_iP=\begin{pmatrix}C_i^c \!\!\!& \ast\end{pmatrix}\!,
$$
where $A_i^c$ and $B_i^c$ belong to $M_r(\mathbb{R})$ and $M_{r,m}(\mathbb{R})$ respectively. Moreover, the switched linear control system defined on $\mathbb{R}^r$ associated with $\mathcal{S}^{\rm arb}(\mathscr{M}^c)$, where 
$$
\mathscr{M}^c=\{(A_1^c,B_1^c,C_1^c),\dots,(A_k^c,B_k^c,C_k^c)\},
$$
is $\mathcal{F}^{\rm arb}(\mathscr{M}^c)$-controllable. We refer to $\mathscr{M}^c$ as the \emph{controllable part} of 
$\mathscr{M}$. Notice that  
$\mathscr{M}^c=\Pi^c(\mathscr{M})$ where
%\begin{equation}\label{eq:pi-c}
\[
\Pi^c(A,B,C)=(UP^{-1}APU^T,UP^{-1}B,CPU^T),
\]
%\end{equation}
with $U=(\Id_r\ 0_{r,n-r})$.
Also notice that the output  $y$  corresponding to the original system is equal to $y=C^cx^c$ and thus the original switched  linear control system has the same $L_2$-gain as the one reduced to the reachable space. 

Similarly, there exists a direct sum decomposition of the state space 
$\mathbb{R}^n=\mathfrak{O}(\mathcal{G})\oplus F$ and  an invertible $n\times n$ matrix $Q$ such that, if $s=\dim \mathfrak{O}(\mathcal{G})$ and $Q^{-1}x=(x^o,x^u)$ with 
$x^o\in\mathfrak{O}(\mathcal{G})$ and $x^u\in F$, one has for $1\leq i\leq k$,
\begin{align*}
Q^{-1}\!A_iQ&=\begin{pmatrix}A_i^o & 0\\ \ast & A_i^u\end{pmatrix}\!,\ 
Q^{-1}\!B_i=\begin{pmatrix}B_i^o\\ \ast\end{pmatrix}\!,\ C_iQ=\begin{pmatrix}C_i^o\!\!& \! 0\end{pmatrix}\!,
\end{align*}
where  $A_i^o$ and $C_i^o$ belong to $M_s(\mathbb{R})$ and $M_{p,s}(\mathbb{R})$ respectively. Moreover, the switched linear control system defined on $\mathbb{R}^s$ associated with $\mathcal{S}^{\rm arb}(\mathscr{M}^o)$, where 
$$
\mathscr{M}^o=\{(A_1^o,B_1^o,C_1^o),\dots,(A_k^o,B_k^o,C_k^o)\},
$$
is $\mathcal{F}^{\rm arb}(\mathscr{M}^o)$-observable. We refer to $\mathscr{M}^o$ as the \emph{observable part} of 
$\mathscr{M}$. Notice that  
$\mathscr{M}^o=\Pi^o(\mathscr{M})$ where
%\begin{equation}\label{eq:pi-o}
\[
\Pi^o(A,B,C)=(VQ^{-1}AQV^T,VQ^{-1}B,CQV^T),
\]
%\end{equation}
and $V=(\Id_s\ 0_{s,n-s})$.
The corresponding output  $y$ being equal to $y=C^ox^o$, one deduces the equality of the $L_2$-gains of the original switched  linear control system and of the one reduced to the observable space.

From the above, one can proceed as follows. Consider a switched linear control system \eqref{LSCS-D} associated with  $\mathcal{S}^{\rm arb}(\mathscr{M})$, where 
$\mathscr{M}=\{(A_1,B_1,C_1),\dots,(A_k,B_k,C_k)\}$. One first reduces it to its reachable  space $\mathfrak{R}(\mathcal{F}^{\rm arb}(\mathscr{M}))$. 
We thus get a $\mathcal{F}^{\rm arb}(\Pi^c(\mathscr{M}))$-controllable switched linear control system  with same 
$L_2$-gain as that of the original system. Then, one reduces the latter system to its observable space
$\mathfrak{O}(\mathcal{F}^{\rm arb}(\Pi^c(\mathscr{M})))$ to finally obtain a switched linear control system associated with $\mathcal{S}^{\rm arb}(\mathscr{M}^{\rm min})$ 
where $\mathscr{M}^{\rm min}=\Pi^o(\Pi^c(\mathscr{M}))$.
The latter system is finally $\mathcal{F}^{\rm arb}(\mathscr{M}^{\rm min})$-controllable and $\mathcal{F}^{\rm arb}(\mathscr{M}^{\rm min})$-observable and 
\begin{equation}\label{L2-real}
\gamma_2(\mathcal{S}^{\rm arb}(\mathscr{M}^{\rm min}))=\gamma_2(\mathcal{S}^{\rm arb}(\mathscr{M})).
\end{equation}
We refer to $\mathcal{S}^{\rm arb}(\mathscr{M}^{\rm min})$ and the corresponding switched linear control system as a \emph{minimal realization} for the linear switched linear control system associated with $\mathcal{S}^{\rm arb}(\mathscr{M})$ and \eqref{LSCS-D}. We also say that $n'$ is the \emph{dimension} of such a minimal realization. 

\begin{remark}\label{rem:realization}
Note that even though the dimension $n'$ of $\mathfrak{O}(\mathcal{F}^{\rm arb}(\Pi^c(\mathscr{M})))$ is uniquely defined by the original switched linear control system, the minimal realization is not unique since it depends on the choice of supplementary spaces to  $\mathfrak{R}(\mathcal{F}^{\rm arb}(\mathscr{M}))$ in $\R^n$ and to $\mathfrak{O}(\mathcal{F}^{\rm arb}(\Pi^c(\mathscr{M})))$ in $\mathfrak{R}(\mathcal{F}^{\rm arb}(\mathscr{M}))$. However, one deduces from Eq.~\eqref{L2-real} that 
$\gamma_2(\mathcal{S}^{\rm arb}(\mathscr{M}^{\rm min}))$
does not depend on a particular choice of a minimal realization. Moreover, it can be shown that any two minimal realizations with switching signal value sets $\mathscr{M}^{\rm min}_1$ and $\mathscr{M}^{\rm min}_2$ are algebraically similar, i.e., there exists an invertible matrix $G\in {\rm GL}_{n'}(\R)$ so that 
\begin{equation}\label{eq:realization}
\mathscr{M}^{\rm min}_2=\{(G^{-1}AG,G^{-1}B,CG)\mid (A,B,C)\in \mathscr{M}^{\rm min}_1\}.
\end{equation}
All the results presented in this paragraph belong to the theme of realization theory of switched linear control systems and we refer to \cite{petreczky} for a thorough presentation of such a theory.  
\end{remark}

Finally, it must be recalled that \cite{Sun-Ge-Lee} also provides a nice and explicit geometric description of $\mathfrak{R}(\mathcal{F}^{\rm arb}(\mathscr{M}))$ and $\mathfrak{O}(\mathcal{F}^{\rm arb}(\mathscr{M}))$ in terms of the data of the problem. 
Let us recall here the details of such a geometric description for $\mathfrak{R}(\mathcal{F}^{\rm arb}(\mathscr{M}))$ (the corresponding results for  $\mathfrak{O}(\mathcal{F}^{\rm arb}(\mathscr{M}))$ being standardly derived by  duality). 

We first need the following notation. If $A$ is an $n\times n$ matrix and $\mathcal{B}$ is a subspace of $\mathbb{R}^n$, let $\Gamma_A\mathcal{B}$ be the subspace of $\mathbb{R}^n$ given by
$$
\Gamma_A\mathcal{B}=\mathcal{B}+A\mathcal{B}+\dots+A^{n-1}\mathcal{B}.
$$
For $1\leq i\leq k$, let $\mathcal{D}_j=\hbox{Im }[B_i\ A_iB_i\ \dots\ A_i^{n-1}B_i]$. Moreover (see \cite[Section 3.1]{Sun-Ge-Lee}), define 
recursively the sequence of subspaces of $\mathbb{R}^n$ denoted $\mathcal{V}_j$, $j\geq 1$, by 
\begin{eqnarray*}
\mathcal{V}_1&=&\mathcal{D}_1+\dots+\mathcal{D}_k,\\
\mathcal{V}_{j+1}&=&\Gamma_{A_1}\mathcal{V}_j+\dots+\Gamma_{A_k}\mathcal{V}_j,
\end{eqnarray*}
and finally set $\mathcal{V}(\mathscr{M})=\sum_{j\geq 1}\mathcal{V}_j$. From the variation of constants formula, it is not difficult to see that $\mathfrak{R}(\mathcal{F}^{\rm arb}(\mathscr{M}))$ is included in $\mathcal{V}(\mathscr{M})$. The converse inclusion is also true but more delicate to establish, cf.~\cite[Theorem 1]{Sun-Ge-Lee} and the proof of Proposition~\ref{prop-reach} below.
In the subsequent paragraphs, we will use these results to derive the existence of a minimal realization associated with \eqref{LSCS-D} and any class $\mathcal{T}$ 
considered in Section~\ref{s:classes} together with its corresponding concatenation of subfamilies of switching signals introduced in Section~\ref{f:conca}. 

We now generalize the above construction to a bounded set 
$\mathscr{M}\subset M_n(\mathbb{R})\times  M_{n, m}(\mathbb{R})\times  M_{p,n}(\mathbb{R})$. We associate a subspace $\mathcal{V}(\mathscr{M})$ of  $\mathbb{R}^n$ as follows. First, consider 
\begin{align*}
\mathcal{V}_1(\mathscr{M})=\hbox{Span} \{A^jb_l\mid & \ 0\leq j\leq n-1,\ 1\leq l\leq m,\\ 
&(A,[b_1\ldots b_m])\in \pi_{A\times B}(\mathscr{M})\}.
\end{align*}
Then, define recursively for $j\geq 1$
\begin{align*}
\mathcal{V}_{j+1}(\mathscr{M})=\hbox{Span} \{A^jv\mid & \ 0\leq j\leq n-1,\\ 
&A\in\pi_A(\mathscr{M}),\ v\in \mathcal{V}_j(\mathscr{M})\},
\end{align*}
and finally set $\mathcal{V}(\mathscr{M})=\sum_{j\geq 1}\mathcal{V}_j(\mathscr{M})$. Taking a (finite) generating family of $\mathcal{V}(\mathscr{M})$, one can extract a finite subset 
$\mathscr{M}^{\rm finite}$
of $\mathscr{M}$ such that $\mathcal{V}(\mathscr{M}^{\rm finite})=\mathcal{V}(\mathscr{M})$. Hence
\begin{align}\label{eq:equality}
\mathcal{V}(\mathscr{M})&\supset \mathfrak{R}(\mathcal{F}^{\rm arb}(\mathscr{M}))\supset\mathfrak{R}(\mathcal{F}^{\rm arb}(\mathscr{M}^{\rm finite}))\nonumber\\
&=\mathcal{V}(\mathscr{M}^{\rm finite})=\mathcal{V}(\mathscr{M}),
 \end{align}
 where the first inclusion is deduced from the variation of constants formula. Therefore, all the sets appearing in \eqref{eq:equality} coincide.

We thus prove the following proposition.
\begin{proposition}\label{prop-reach}
Consider a switched linear control system of the type \eqref{LSCS-D} associated with a class $\mathcal{T}$ of measurable switching laws taking values in a bounded set $\mathscr{M}\subset M_n(\mathbb{R})\times  M_{n, m}(\mathbb{R})\times  M_{p,n}(\mathbb{R})$.  Let $(\mathcal{T},\mathcal{G})$ satisfy Assumption 
{\bf T2}. Then 
$\mathfrak{R}(\mathcal{G})=\mathfrak{R}(\mathcal{F}^{\rm arb}(\mathscr{M}))=\mathcal{V}(\mathscr{M})$ and there exist $t^*>0$ and a switching law
$(A,B,C)\in \mathcal{G}_{t^*}$ such that the range of $W^c(t^*)$, the controllability Gramian in time $t^*$ associated with $\dot x(t) =A(t)x(t)+B(t)u(t)$, is equal to $\mathcal{V}(\mathscr{M})$.
\end{proposition}
\begin{IEEEproof} First notice that the equality $\mathfrak{R}(\mathcal{F}^{\rm arb}(\mathscr{M}))=\mathcal{V}(\mathscr{M})$ is contained in Eq.\eqref{eq:equality} and one has the trivial inclusion $\mathfrak{R}(\mathcal{G})\subset\mathcal{V}(\mathscr{M})$. 
Let us prove the opposite inclusion.

As done above there exists a finite subset 
$\mathscr{M}^{\rm finite}$
of $\mathscr{M}$ such that $\mathcal{V}(\mathscr{M}^{\rm finite})=\mathcal{V}(\mathscr{M})$. It is proved in \cite[Theorem 1]{Sun-Ge-Lee} that there exists a piecewise-constant periodic switching law
$(A,B,C)$ taking values in 
$\mathscr{M}^{\rm finite}$ and a time $t^*>0$ such that  the range of $W^c(t^*)$, the controllability Gramian in time $t^*>0$ associated with $\dot x(t) =A(t)x(t)+B(t)u(t)$, is equal to $\mathcal{V}(\mathscr{M})$. Fix $\varepsilon>0$ such that if  an $n\times n$ matrix $W$ satisfies $\Vert W-W^c(t^*)\Vert< \varepsilon$ then the rank of $W$ is larger than or equal to the rank of $W^c(t^*)$. Let $\lambda\mapsto (A_\lambda,B_\lambda,C_\lambda)$ be the path provided by Assumption {\bf T2}. Then there exists 
$\lambda^*\in (0,1]$ such that $(A_{\lambda^*},B_{\lambda^*},C_{\lambda^*})\in \mathcal{G}$ and
the rank of the corresponding controllability Gramian $W_{\lambda^*}^c(t_{\lambda^*})$ is larger than or equal to the rank of $W^c(t^*)$, itself equal to $\dim\mathcal{V}(\mathscr{M})$. Since the range of $W_{\lambda^*}^c(t_{\lambda^*})$ is included in $\mathfrak{R}(\mathcal{G})$ this concludes the proof of the proposition.
\end{IEEEproof}

An analoguous statement is obtained regarding observability spaces and, following the first part of the section, one finally derives the subsequent result.
\begin{proposition}\label{prop:reduction}
Consider a switched linear control system of the type \eqref{LSCS-D} associated with a class $\mathcal{T}$ of measurable switching laws taking values in a bounded set $\mathscr{M}\subset M_n(\mathbb{R})\times  M_{n, m}(\mathbb{R})\times  M_{p,n}(\mathbb{R})$.  Let $(\mathcal{T},\mathcal{G})$ satisfy Assumption 
{\bf T2}.
Let $n'$ be the dimension of any minimal realization of the switched linear control system associated with $\mathcal{S}^{\rm arb}(\mathscr{M})$.
Pick one such minimal realization and consider the corresponding 
surjective linear mapping $\Pi$ from $M_n(\mathbb{R})\times  M_{n, m}(\mathbb{R})\times  M_{p,n}(\mathbb{R})$ to $M_{n'}(\mathbb{R})\times  M_{n', m}(\mathbb{R})\times  M_{p,n'}(\mathbb{R})$
%$$
%\Pi:M_n(\mathbb{R})\times  M_{n, m}(\mathbb{R})\times  M_{p,n}(\mathbb{R})\rightarrow M_{n'}(\mathbb{R})\times  M_{n', m}(\mathbb{R})\times  M_{p,n'}(\mathbb{R})
%$$
describing its matrix representation.
Denote by $\mathcal{T}^{\rm min}$ the class  
$$
\mathcal{T}^{\rm min}=\{t\mapsto\Pi(A(t),B(t),C(t))\mid  (A,B,C)\in \mathcal{T}\}
$$
taking values in $\mathscr{M}^{\rm min}=\Pi(\mathscr{M})$. 
Then the 
switched linear control system corresponding to $\mathcal{T}^{\rm min}$ 
is $\mathcal{G}^{\rm min}$-controllable and $\mathcal{G}^{\rm min}$-observable in the sense of 
Definition~\ref{rea-obs} with $\mathcal{G}^{\rm min}=\Pi(\mathcal{G})$. 
Moreover, $(\mathcal{T}^{\rm min},\mathcal{G}^{\rm min})$ satisfies Assumption 
{\bf T2} and $\gamma_2(\mathcal{T})=\gamma_2(\mathcal{T}^{\rm min})$.
Finally, if $(\mathcal{T},\mathcal{G})$ satisfies Assumption~{\bf T1} then $(\mathcal{T}^{\rm min},\mathcal{G}^{\rm min})$ does.
\end{proposition}
\begin{IEEEproof}
This proposition is a simple consequence of the construction of minimal realizations in the case of arbitrary switching, as described in the first part of the section, and of Proposition~\ref{prop-reach}. 
\end{IEEEproof}

For $\mathcal{T}$ and $\mathcal{T}^{\rm min}$ as in the statement of Proposition~\ref{prop:reduction}, we say that 
the 
switched linear control system corresponding to $\mathcal{T}^{\rm min}$ 
is a \emph{minimal realization} of the switched linear control system associated with $\mathcal{T}$.
It follows from Remark~\ref{rem:realization} that any two such minimal realizations 
are algebraically similar.

\subsection{Finiteness of the $L_2$-gain}
Consider the switched linear control system \eqref{LSCS-D} 
and the corresponding class  $\mathcal{T}$ of switching laws with values in the bounded set $\mathscr{M}$ 
of triples of matrices. 
Let us introduce the following definition.
\begin{Definition}\label{def:UOC}
We say that \eqref{LSCS-D}  is \emph{uniformly observable} (or, equivalently, that $\mathcal{T}$ is \emph{uniformly observable}) if there exist $T,\gamma>0$ such that, for every 
$(A,B,C)\in \mathcal{T}$ and every $t\geq 0$, one has $W^o(t,t+T)\geq \gamma \Id_n$, where $W^o(t,t+T)$ is the observability Gramian in time $T$ associated with $(A(t+\cdot),B(t+\cdot),C(t+\cdot))$. 
\end{Definition}
The following remark will be used in the sequel.

\begin{remark}\label{closure} 
If  \eqref{LSCS-D}  is uniformly observable and $T,\gamma$ are as in Definition~\ref{def:UOC}, 
 the observability Gramian
$W^o(t,t+T)$ in time $T$ associated with a switching law belonging to the weak-$\star$ closure of  $\mathcal{T}$ still satisfies $W^o(t,t+T)\geq \gamma \Id_n$. 
 \end{remark}

\begin{remark}\label{compact} 
Consider the case where $\mathcal{T}$ contains  all the constant $\mathscr{M}$-valued switching signals. 
(Notice that this is the case for all classes of switching signals introduced in Section~\ref{s:classes} except that of persistently exciting signals.) It is then easy to see that uniform observability implies that $(A,C)$ is observable for every $(A,B,C)\in \mathscr{M}$.
 In the case where $\mathscr{M}$ is compact and the class of signals under consideration is $\mathcal{S}^{{\rm d},\tau}(\mathscr{M})$ 
for some $\tau>0$, one easily shows that the converse is also true, namely, the observability of each pair $(A,C)$  implies uniform observability.
Let us stress that the uniform observability assumption is weaker than the minimal realization 
assumption (every triple $(A,B,C)\in \mathscr{M}$ is a minimal realization, i.e., $(A,B)$ is controllable and $(A,C)$ is observable) needed in \cite{Hespanha-IEEE,Margaliot-Hespanha}.
 \end{remark}
 
We can now state the main result of this section.
\bt\label{th1}
Assume that \eqref{LSCS-D} 
admits a minimal realization defined on $\R^{n'}$, $n'\leq n$, which is given by $\dot x^{\rm min}(t)=A^{\rm min}(t)x(t)+B^{\rm min}(t)u(t)$ with output $y^{\rm min}(t)=C^{\rm min}(t)x^{\rm min}(t)$  and associated with a class $\mathcal{T}^{\rm min}$ of switching signals taking values in $\mathscr{M}^{\rm min}$, and a family $\mathcal{G}^{\rm min}$ satisfying the following assumptions:
\begin{itemize} %\begin{description} 
\item[$(a)$] $(\mathcal{T}^{\rm min},\mathcal{G}^{\rm min})$ 
satisfies Assumptions~{\bf T1} and {\bf T2};
\item[$(b)$] for every subspace $V$ of $\R^{n'}$ and every
linear projection $P_V:\R^{n'}\rightarrow V$, the class of signals  $P_V\mathcal{G}_A^{\rm min}P_V^{\#}$ satisfies Assumption {\bf A3}, where 
 we use $P_V^{\#}$ to denote the dual map from $V$ to $\R^n$ defined by $x^TP_V y=y^TP^{\#}_Vx$ for every $x\in V$ and $y\in\R^n$.
\end{itemize} %\end{description} 
Then $\gamma_2(\mathcal{T})$ is finite if $\rho(\mathcal{T}_A^{\rm min})<1$ and infinite if either $\rho(\mathcal{T}_A^{\rm min})>1$ or $\rho(\mathcal{T}_A^{\rm min})=1$ and 
$\mathcal{T}^{\rm min}$ is uniformly observable.
\et
\begin{IEEEproof} Thanks to Proposition~\ref{prop:reduction}, it is enough to treat the case 
$\mathcal{T}=\mathcal{T}^{\rm min}$.

Assume first that $\rho(\mathcal{T}_A)<1$. Taking into account the definition  of $\rho(\mathcal{T}_A)$ (see \eqref{eq:bohl}) and the boundedness of $\pi_A(\mathscr{M})$, one gets the following exponential decay estimate: 
 there exist $K_1,\lambda>0$ such that, for every $A\in \mathcal{T}_A$ and every $0\leq s\leq t$, one has
$$
\left\|\Phi_{A}(t,s)\right\|\leq K_1e^{-\lambda(t-s)}.
$$
As a consequence of the above and the boundedness of $\mathscr{M}$, one deduces that 
there exists $K_2>0$ such that, for every $u\in L_2$, $(A,B,C)\in \mathcal{T}$ and 
$t\geq 0$, one has
\begin{equation}\label{eq:ics}
\Vert y_u(t)\Vert\leq K_2\int_0^t e^{-\lambda(t-s)}\Vert u(s)\Vert ds.
\end{equation}
If $\chi_{[0,+\infty)}$ denotes the characteristic function of $[0,+\infty)$, the integral function on the right-hand side of \eqref{eq:ics} can be interpreted as the convolution
of 
$$
f_1(t)=\chi_{[0,+\infty)}(t)e^{-\lambda t},\quad f_2(t)=\chi_{[0,+\infty)}(t) \Vert u(t)\Vert. 
$$
That yields at once that $\|y_u\|_2\leq \frac{K_2}{\lambda} \|u\|_2$, hence the conclusion.

Assume now that 
$\rho(\mathcal{T}_A)\geq 1$. 
It is well-known
(see, e.g., \cite[Proposition 2]{CMS1} for details), 
that, up to a common linear change of coordinates,
 every matrix in $A\in \mathscr{M}_A$ has the  upper triangular block form 
%\begin{equation}
\[
A=\left(\ba{ccccc} A_{11} & A_{12} & \cdots & &  \\ 0 & A_{22} & A_{23} & \cdots  & 
\\ 0 & 0 & A_{33} & A_{34} & \cdots  \\ \vdots & \ddots & \ddots &  \ddots  & \\ 0 & 
\cdots & \cdots & 0  &  A_{qq} \ea\right)\,,
%\label{block}
\]
%\end{equation}
where, for $i=1,\ldots,q$, each $A_{ii}$ is in $M_{n_i-n_{i-1}}(\mathbb{R})$, $n_i\in\mathbb{N}$ and the set $\A_i:=\{A_{ii}\mid A\in\mathscr{M}_A\}$ is irreducible (whenever  $\mathscr{M}_A$ is irreducible one has $q=1$ and $\A_1=\mathscr{M}_A$). Define 
the \emph{subsystems of $\mathscr{M}_A$} 
as the switched systems corresponding to the  sets $\A_i$ and the class of switching signals $\mathcal{T}_{A,i}:=\{A_{ii}\mid A\in\mathcal{T}_A\}$ and $\mathcal{G}_{A,i}:=\{A_{ii}\mid A\in\mathcal{G}_A\}$. One can then show that 
 $(\mathcal{T}_{A,1},\mathcal{G}_{A,1}),\dots,(\mathcal{T}_{A,q},\mathcal{G}_{A,q})$ satisfy Assumptions~{\bf A0--A3}. Indeed, Assumptions {\bf A0} and {\bf A1} follow from Assumption {\bf T1}; 
Assumption {\bf A2} is a consequence of Assumption {\bf T2}  and Proposition~\ref{prop:A2}
while Item $(b)$ yields Assumption {\bf A3}.

Moreover, an induction argument on the number of subsystems and a standard use of the variation of constants formula yields 
$$
\rho(\mathcal{T}_A)=\max_{1\leq i\leq q}\rho(\mathcal{T}_{A,i}).
$$
Let ${\bar \imath}\leq q$ be the largest index such that $\rho(\mathcal{T}_{A,{\bar \imath}})=\rho(\mathcal{T}_A)$. 
Since \eqref{LSCS-D} 
is $\mathcal{G}$-controllable, 
there exists a time $\bar{t}> 0$, a switching law $(A^0,B^0,C^0)\in\mathcal{G}_{\bar{t}}$ and a measurable control $\bar{u}$ defined on $[0,\bar{t}]$ so that the corresponding trajectory $x_{\bar{u}}$ starting at $0$ reaches some point $\bar{x}=(0,\ldots,0,\bar{x}_{\bar \imath},0,\ldots,0)^T\ne 0$ in time $\bar{t}$.

Since $(\mathcal{T}_{A,{\bar\imath}},\mathcal{G}_{A,{\bar \imath}})$ satisfies the hypotheses of Theorem~\ref{anticipazioni},  the corresponding semigroup $\Mexp(\mathcal{G}_{A,{\bar\imath}})$ is quasi-Barabanov and admits a quasi-extremal trajectory $x^{\rm q-ext}_{\bar \imath}$
starting at $\bar{x}_{\bar \imath}$ and corresponding to a signal $A_{{\bar \imath}{\bar \imath}}$ belonging to the weak-$\star$ closure of $(\mathcal{G}_{A,{\bar \imath}})_\infty$. 
Let $\bar{x}$ be a trajectory of $\dot x=A(t)x(t)$ starting at $\bar{x}$ and corresponding to some signal $\bar{A}$ in the closure of $({\mathcal{G}}_A)_\infty$ so that
$\bar{A}_{{\bar \imath}{\bar \imath}}=A_{{\bar \imath}{\bar \imath}}$.
Notice that $\bar{x}_{\bar \imath}=x^{\rm q-ext}_{\bar \imath}$.

Let $\widebar{E}$ be  the signal defined as the concatenation of 
$A^0$ and $\bar{A}$ and let $A^0_{{\bar \imath}{\bar \imath}}$ and $\widebar{E}_{{\bar \imath}{\bar \imath}}$ be the $({\bar \imath},{\bar \imath})$-blocks corresponding to $A^0$  and $\widebar{E}$ respectively. 
Let $(E^l)_{l\geq 0}$ (and then $(E^l_{{\bar \imath}{\bar \imath}})_{l\geq 0}$ ) be a sequence in $({\mathcal{G}}_A)_\infty$ ($({\mathcal{G}}_{A,{\bar \imath}})_\infty$ respectively) weak-$\star$ converging to $\bar{A}$ ($\bar{A}_{{\bar \imath}{\bar \imath}}$ respectively). 
Let $(F^l,G^l)_{l\geq 0}$ be such that $(E^l,F^l,G^l)\in \mathcal{G}_\infty$ for every $l\geq 0$. 
The sequence $(A^0\ast E^l)_{l\geq 0}$ ($(A^0_{{\bar \imath}{\bar \imath}}\ast E^l_{{\bar \imath}{\bar \imath}})_{l\geq 0}$) is in  $({\mathcal{G}}_A)_\infty$ ($({\mathcal{G}}_{A,{\bar \imath}})_\infty$ respectively) and 
weak-$\star$ converges to $\widebar{E}$ ($\widebar{E}_{{\bar \imath}{\bar \imath}}$ respectively). Moreover, for every positive time $t\ge \bar t$, the convergence of $s\mapsto \Phi_{E_l}(s-\bar t,0)$ to $s\mapsto \Phi_{\widebar E}(s,\bar t)$ is uniform with respect to $s\in [\bar t,t]$ and similarly for the corresponding $({\bar \imath},{\bar \imath})$-blocks.

Denote by $x^l$ the trajectory of \eqref{LSCS-D} starting at the origin and associated with the control equal to 
$\bar{u}$ on $[0,\bar t]$ and zero on $[\bar t,+\infty)$, 
and with the switching signal $(A^0\ast E^l,B^0\ast F^l,C^0\ast G^l)$.  Let $y^l$ be the corresponding output.
Notice that $x^l(s)=\Phi_{E_l}(s-\bar t,0) \bar x$ for every $s\geq \bar t$. 
For $l$ large enough and for every $s\in [\bar t,t]$, one has
%\begin{equation}
\begin{align*}%\label{est-L2}
\Vert x^l(s)\Vert%=\Vert \cexp_{\bar t}^{s} E^l(\tau)d\tau \bar x\Vert 
&\geq
\Vert{\Phi_{E^l_{{\bar \imath}{\bar \imath}}}(s-\bar t,0)} \bar{x}_{\bar \imath}\Vert\nonumber\\
&\geq \frac{1}2\Vert\Phi_{\widebar{E}_{{\bar \imath}{\bar \imath}}} (s,\bar{t})\bar{x}_{\bar \imath}\Vert
\geq \frac{1}{2C_{\rm qx}}
\rho(\mathcal{T}_A)^{s-\bar t}\Vert\bar{x}_{\bar \imath}\Vert.
\end{align*}
%\end{equation}
 Hence, there exists a positive constant $K_3$ independent of $t\geq \bar t$ such that, for every 
$l$ large enough and every
$s\in [\bar t,t]$, one has 
\begin{equation}\label{est-L22}
\Vert x^l(s)\Vert\geq K_3 \rho(\mathcal{T}_A)^{s-\bar t}.
\end{equation}

Assume now that $\rho(\mathcal{T}_A)=1$ and 
$\mathcal{T}$ is uniformly observable, and let $T,\gamma$ be as in Definition~\ref{def:UOC}. 
Given two positive integers $l,j$, let
$W^o_{l,j}$ be the observability Gramian in time $T$ associated with $(E^l(jT+\cdot),F^l(jT+\cdot),G^l(jT+\cdot))$.
Applying the uniform observability assumption to each $(E^l(jT+\cdot),F^l(jT+\cdot),G^l(jT+\cdot))$ and taking into account Eq.~\eqref{est-L22}, 
we get, for $J\in \N$ and $l$ large enough,
\begin{align*}
\gamma_2(\mathcal{T})&\ge \frac{\| y^l\|_2^2}{\Vert \bar{u}\Vert_2^2}\geq   \frac{\int_0^{J T} 
\Vert G^l(t)x^l(\bar t+t)\Vert^2dt }{\Vert \bar{u}\Vert_2^2}\\
&= 
\frac{\sum_{j=0}^{J-1}\int_{jT}^{(j+1)T}
\Vert G^l(t)x^l(\bar t+t)\Vert^2dt }{\Vert \bar{u}\Vert_2^2}\\
&= 
\frac{\sum_{j=0}^{J-1}
x^l(\bar t+jT)^TW^o_{l,j}x^l(\bar t+jT)}{\Vert \bar{u}\Vert_2^2}\\
&\geq \gamma\frac{\sum_{j=0}^{J-1}
\Vert x^l(\bar t+jT)\Vert^2 }{\Vert \bar{u}\Vert_2^2}\geq \gamma J\frac{K_3^2}{\Vert \bar{u}\Vert_2^2}.
\end{align*}
This implies that the $L_2$-gain $\gamma_2(\mathcal{T})$ is infinite. 

Let now $\rho(\mathcal{T}_A)>1$.
Let $\bar\imath$, $\bar u$ and $\bar x$ be as above. 
Since the semigroup $\Mexp(\mathcal{G}_{A,{\bar\imath}})$ is quasi-Barabanov,  
we can find for every $l\in \N$ a time $t_l\geq l$ and a switching law 
$(\hat E^l,\hat F^l,\hat G^l)\in \mathcal{G}_{t_l}$  such that 
\begin{equation}\label{est-L23}
  \Vert\Phi_{\hat E^l}(t_l,0) \bar{x}
\Vert\geq 
C_{\rm qb}\rho(\mathcal{T}_A)^{t_l}\Vert\bar{x}
\Vert\geq C_{\rm qb}\rho(\mathcal{T}_A)^l\Vert\bar{x}
\Vert.
\end{equation}
According to {\bf T2} and because of the observability counterpart of Proposition~\ref{prop-reach}, there exist $\tilde s>0$
and a switching law $(A_*,B_*,C_*)\in \mathcal{G}_{\tilde s}$ such that 
the observability Gramian  $W^o(\tilde s)$ in time $\tilde s$ associated with $\dot x(t) =A_*(t)x(t)+B_*(t)u(t)$, $y(t)=C_*(t)x(t)$, is invertible.
Up to a suitable extension on $(\tilde s,\infty)$, we can assume that  $(A_*,B_*,C_*)$ belongs to $\mathcal{G}_{\infty}$. 
For $l\geq 0$, consider the sequence of switching signals 
$$
S^l=(A^0\ast\hat  E^l\ast A_*,B^0\ast \hat F^l\ast B_*,C^0\ast \hat G^l\ast C_*)\in \mathcal{G}_{\infty}. 
$$
Denote by $x^l$ the corresponding 
trajectory of Eq.~\eqref{LSCS-D} starting at the origin and associated with the control equal to $\bar{u}$ on $[0,\bar t]$ and zero for $t> \bar t$ and let $y^l$ be the corresponding
 output. 
Note that $x^l(\bar t+t_l)=\Phi_{\hat E^l}(t_l,0) \bar{x}$. 
It then follows from \eqref{est-L23} that there exists a positive constant $K_4$ independent of $l$ so that 
\begin{align*}
\gamma_2(\mathcal{T})&\ge \frac{\| y^l\|_{2,t_l}^2}{\Vert \bar{u}\Vert_2^2}\geq \frac{\int_{\bar t+t^l}^{\bar t+t^l+\tilde s}\Vert C_*(s)x^l(s)\Vert^2ds}{\Vert \bar{u}\Vert_2^2}\\
&= \frac{x^l(\bar t+t_l)^TW^o(\tilde s)x^l(\bar t+t_l)}{\Vert \bar{u}\Vert_2^2}\geq K_4  \rho(\mathcal{T}_A)^{2l},
\end{align*}
and the right-hand side clearly tends to infinity as $l$ tends to infinity.
\end{IEEEproof}

\begin{remark}
Note that the  value of $\rho(\mathcal{T}_A^{\rm min})$ does not depend on the particular choice of the minimal realization thanks to Eq.~\eqref{eq:realization}.
\end{remark} 

 Under the hypotheses of Proposition~\ref{lem:A4}, the classes of switching signals considered in Section~\ref{s:classes} together with the corresponding families of Section~\ref{f:conca} satisfy all the hypotheses of Theorem~\ref{th1}. 
 As a consequence, we can now answer some of the questions raised by Hespanha in \cite{Hespanha-unsolved}.
\bt\label{q-hespanha}
Let  $\mathscr{M}$ be a bounded subset of $M_n(\R)\times M_{n,m}(\R)  \times  M_{p,n}(\mathbb{R})$ with $n,m,p$ positive integers and $\tau\geq 0$. Consider the switched linear control system $\dot x(t)=A(t)x(t)+B(t)u(t)$, $y(t)=C(t)x(t)$, where the switching signal $(A,B,C)$ belongs to the class  $\mathcal{S}^{{\rm d},\tau}(\mathscr{M})$ of piecewise constant signals {with} \emph{dwell-time} $\tau$. Let $\gamma_2(\tau)$ be the $L_2$-gain associated with 
 $\mathcal{S}^{{\rm d},\tau}(\mathscr{M})$. Consider a minimal realization defined on $\R^{n'}$, $n'\leq n$, given by $\dot x^{\rm min}(t)=A^{\rm min}(t)x(t)+B^{\rm min}(t)u(t)$ with output $y^{\rm min}(t)=C^{\rm min}(t)x^{\rm min}(t)$ and associated with a class $\mathcal{S}^{{\rm d},\tau}(\mathscr{M}^{\rm min})$ where $\mathscr{M}^{\rm min}$ is a bounded subset of $M_{n'}(\R)\times M_{n',m}(\R) \times  M_{p,n'}(\mathbb{R})$. Assume furthermore that this minimal realization is uniformly observable.
  
Then, $\gamma_2(\tau)$ is finite if and only if $\rho(\mathcal{S}^{{\rm d},\tau}(\mathscr{M}_A^{\rm min}))<1$ and, if $ \tau_{\min}$ is defined as
 \begin{align*}
 \tau_{\min}=  \inf\{\tau>0\mid \gamma_2(\tau)\mbox{ is finite}\},
 \end{align*}
one has the following characterization: $ \tau_{\min}=\inf\{\tau>0\mid \rho(\mathcal{S}^{{\rm d},\tau}(\mathscr{M}_A^{\rm min}))<1\}$. 
\et

\begin{remark}
The %statement of the 
%conclusion of the 
theorem
still holds true if we replace the uniformly observability assumption   by 
the hypothesis that there exists at most one $\tau>0$ such that 
$\rho(\mathcal{S}^{{\rm d},\tau}(\mathscr{M}_A^{\rm min}))=1$. 
\end{remark}

\begin{remark}
One can derive results similar to the previous theorem when one considers variations of certain parameters used in the definition of classes other than $\mathcal{S}^{{\rm d},\tau}(\mathscr{M})$. For instance, one can characterize the set of $\mu\in (0,T)$ such that $\gamma_2(\mathcal{S}^{{\rm pe},T,\mu}(\mathscr{M}))$ is finite in terms of the value of $\rho(\mathcal{S}^{{\rm pe},T,\mu}(\mathscr{M}_A^{\rm min}))$.
\end{remark}

\begin{remark}
Recall  that the computation of 
$t_{\min}=\inf\{\tau>0\mid \rho(\mathcal{S}^{{\rm d},\tau}(\mathscr{M}_A^{\rm min}))<1\}$
turns out to be a numerically tractable task. Indeed, \cite{Chesi} proposes an LMI procedure providing a sequence of upper bounds of $t_{\min}$ approximating it arbitrarily well. 
Therefore, 
combining Theorem~\ref{q-hespanha} 
and the LMI-based algorithm of \cite{Chesi} yields 
a numerical procedure for estimating $\tau_{\min}$. 
\end{remark}

Theorem~\ref{th1} shows that,  under the assumption of uniform observability of a minimal realization, the necessary and sufficient condition for finiteness of the $L_2$-gain 
(i.e., generalized spectral radius smaller than one)
is exactly the same as in the unswitched framework. We prove below by means of an example that this is no more the case when the assumption of uniform observability does not hold. For this purpose, we next define a switched linear control system satisfying all the assumptions of Theorem~\ref{th1} (with $\mathcal{T}_{A}=\mathcal{T}_A^{\rm min}=(\mathcal{G}_A^{\rm min})_\infty$) and for which uniform observability does not hold. 
We have the following  example.
\begin{example}\label{prop-partial}
Assume that $\mathcal{T}(\alpha)=\mathcal{S}^{\rm arb}(\mathscr{M}(\alpha))$ where 
$\mathscr{M}(\alpha)=\{(A_i,b_i,c_i)\}_{i=1,2,3}\subset M_3(\R)\times \R^3\times \R^3$, with
\begin{align*}
A_1=&\left(\begin{array}{ccc}  -1 & -\alpha & 0\\   \alpha  & -1 & 0\\ 0 & 0 & -1 \end{array}\right),\ A_2=\left(\begin{array}{ccc} -1 & -\alpha & 0\\   1/\alpha  & -1 & 0\\ 0 & 0 & -1\end{array}\right),\\ 
&A_3=\left(\begin{array}{ccc} -4 & 0 & 1 \\ 0 & -4 & 0 \\1& 0 & -1 \end{array}\right),\\
b_1&=b_2=0,\quad b_3=c_1=c_2=c_3=(0,0,1)^T,
\end{align*}
for $\alpha>0$. We use $\gamma_2(\alpha)$ to denote the $L_2$-gain induced by the switched linear control system given by $\dot x=A(t)x+b(t)u(t)$, $y(t)=c^T(t)x(t)$ and 
$(A,b,c)\in \mathcal{T}(\alpha)$. Then %one has 
we claim that $\mathcal{T}(\alpha)=\mathcal{T}^{\rm min}(\alpha)$ for every $\alpha>0$ 
and
there exists $\alpha_*$ (approximatively equal to $4.5047$) such that $\rho(\mathcal{T}_{A}(\alpha_*))=1$ and $\gamma_2(\alpha_*)\leq 4$.

Before providing a proof, let us note that the assumption of uniform observability does not hold since the observability Gramian in any positive time associated with a switching signal only activating the first two modes contains 
$b_3$ in its kernel.

Let us now prove the claim. Using the results  in~\cite{BBM}, one determines the value $\alpha=\alpha_*\sim 4.5047$ for which the switched system associated with $\mathscr{M}'_A(\alpha)=\{A_1,A_2\}$ is marginally stable (and reducible). In this case, starting from  every point $(x_1,x_2,0)$, there exists  a closed (periodic) $\mathcal{C}^1$ trajectory $\Gamma_{x_1,x_2}$ of the switched system 
lying on the plane $x_3=0$ which can be completely determined by explicit computations. 
In particular, we can pick such a trajectory so that its support
  $\Gamma$ is contained in the set $\{(x_1,x_2,0)\in\mathbb{R}^3\,|\,1\leq x_1^2+x_2^2\leq 3\}$. We define the norm  $v(x_1,x_2,0)$ on the plane $\{x_3=0\}$ by  setting $v^{-1}(1)=\Gamma$. Then $v$  is a Barabanov norm for the restriction of $\mathscr{M}'_A(\alpha_*)$ on  the plane $\{x_3=0\}$. We extend $v$ to a function on $\R^3$, still denoted by $v$,  by setting $v(x_1,x_2,x_3)=v(x_1,x_2,0)$ and it follows by explicit computations that $\|\nabla v(x)\|\leq \sqrt{3}$ and, by homogeneity, that $v(x)=\nabla v(x)^T(x_1,x_2,0)^T$ for every $x\in\R^3$.
Notice, moreover, that
$$v(x)\ge \sqrt{\frac{x_1^2+x_2^2}{3}},\qquad \forall x\in\R^3.$$

Let us consider the positive definite function $V(x)=\frac12 (v(x)^2+x_3^2)$ and observe that $\frac{d }{d t}V(x(t))\leq -x_3(t)^2$ whenever $A(t)=A_1$ or $A(t)=A_2$.
If $A(t)=A_3$, one deduces from the above properties of $v$ that
\begin{align*}
\frac{d }{d t}&V(x(t))=\nabla V(x(t))^T  (A_3 x(t)+u(t)b_3)\\
&= v(x(t)) (\nabla v(x(t))^T A_3 x(t)+u(t)\nabla v(x(t))^Tb_3 )\\
&+ x_3(t) (-x_3(t)+x_1(t)+u(t))\\
&= -4 v(x(t))^2 +v(x(t)) \nabla v(x(t))^T (x_3(t),0,0)^T\\
&+ x_3(t)  (-x_3(t)+x_1(t)+u(t))\\
&\leq  -4v(x(t))^2+2\sqrt{3}v(x(t)) \vert x_3(t)\vert-x_3(t)^2+u(t) x_3(t)\\
&\leq  -x_3(t)^2/4+u(t) x_3(t).
\end{align*}
Hence $\frac{d }{d t}V(x(t))\leq -x_3(t)^2/4+\vert u(t) x_3(t)\vert$ along any trajectory of the switched linear control system.
By integrating the 
previous inequality, using the fact that $x(0)=0$ 
and applying Cauchy--Schwarz inequality, we get 
\begin{align*}
0\leq\liminf_{t\to \infty} V(x(t))&\leq -\frac14 \|x_3\|_2^2+\int_0^{\infty}\vert u(s) x_3(s)\vert ds\\
&\leq -\frac14 \|x_3\|_2^2+\|u\|_2 \|x_3\|_2,
\end{align*}
so that $\|x_3\|_2\leq 4 \|u\|_2$, implying that $\gamma_2(\alpha_*)\leq 4$. This concludes the proof of the claim.
\end{example}

\subsection{Right-continuity and boundedness of the $L_2$-gain}
\label{s:rc}
In this section, we restrict for simplicity our discussion to the class of signals $\mathcal{S}^{{\rm d},\tau}(\mathscr{M})$. (For more general considerations, see Remark~\ref{r-general-r-continuity}.)

Let  $\mathscr{M}$ be a bounded subset of $M_n(\R)\times M_{n,m}(\R)\times  M_{p,n}(\mathbb{R})$ with $n,m,p$ positive integers and $\tau>0$. Consider the switched linear control system $\dot x(t)=A(t)x(t)+B(t)u(t)$, $y(t)=C(t)x(t)$, where the switching law $(A,B,C)$ belongs to the class  $\mathcal{S}^{{\rm d},\tau}(\mathscr{M})$.
 of piecewise constant signals with dwell-time $\tau\ge 0$. 
 Let $\gamma_2(\tau)$ be the $L_2$-gain associated with 
 $\mathcal{S}^{{\rm d},\tau}(\mathscr{M})$.

We can now state our result on the right-continuity of $\tau\mapsto \gamma_2(\tau)$, which 
answers Questions (i) and---partially---(iii) of \cite{Hespanha-unsolved}.

 \begin{proposition}\label{m-prop-cont}
 The function $\gamma_2:[0,\infty)\to [0,\infty]$ is right-continuous, i.e., for every $\bar\tau\in [0,\infty)$, $\lim_{\tau\searrow \bar\tau}\gamma_2(\tau)=\gamma_2(\bar \tau)$.
 \end{proposition}
\begin{IEEEproof}
For every $T>0$ and every $\tau\in[0,\infty)$, define the $L_2$-gain in time $T$ as
\begin{align*}%\label{m-L2G-T}
\gamma_2(\tau,T):=\sup\bigg\{
\frac{\Vert y_{u,\sigma}\Vert_{2,T}}{\Vert u\Vert_{2,T}}\ \vert \ &u\in L_2([0,T],\mathbb{R}^m)\setminus\{0\},\\ 
&\sigma\in \mathcal{S}^{{\rm d},\tau}(\mathscr{M})
\bigg\}.
\end{align*}
It is immediate to see that $\gamma_2(\tau,T)$ is finite for every $(\tau,T)\in[0,\infty)\times (0,\infty)$ and that the maps $\tau\mapsto \gamma_2(\tau,T)$ (for fixed $T>0$) and $T\mapsto \gamma_2(\tau,T)$ (for fixed $\tau\geq 0$) are non-increasing and non-decreasing respectively.
Also notice that $\tau\mapsto \gamma_2(\tau)$ is non-increasing.

We claim that 
\begin{itemize} %\begin{description} 
\item[(i)] $\lim_{T\rightarrow\infty}\gamma_2(\tau,T)=\gamma_2(\tau)$ for every  $\tau\geq 0$;
\item[(ii)] the map $\tau\mapsto \gamma_2(\tau,T)$ is right-continuous for every $T>0$. 
\end{itemize} %\end{description} 

In order to prove property {\bf (i)} of the claim, notice that, given $\tau\geq 0$, $0<T\leq\infty$, any switching signal $\sigma$ and nonzero control $u\in L_2$, one has 
$$
\lim_{T'\nearrow T}\frac{\Vert y_{u,\sigma}\Vert_{2,T'}}{\Vert u\Vert_{2,T'}}=\frac{\Vert y_{u,\sigma}\Vert_{2,T}}{\Vert u\Vert_{2,T}},
$$
since  $\Vert y_{u,\sigma}\Vert_{2,T'}$ and $\Vert u\Vert_{2,T'}$
converge to $\Vert y_{u,\sigma}\Vert_{2,T}$ and $\Vert u\Vert_{2,T}$ respectively.  
Property {\bf (i)} then follows from 
the definition of $\gamma_2(\tau)$  and the monotonicity of $T\mapsto \gamma_2(\tau,T)$.

Let us now prove point {\bf (ii)} of the claim. 
With $0\le \tau<\tau'$, $T>0$, and $\sigma\in \mathcal{S}^{\rm d,\tau}$ we associate $\sigma'\in L_\infty([0,\infty),\mathscr{M})$  as follows: $\sigma'(\cdot):=\sigma(\xi \cdot)$ where $\xi$ is the largest number in $[0,1]$ such that the restriction of $\sigma'$ to $[0,T]$ has dwell-time $\tau'$. 
Notice that $\xi$ is larger than or equal to $\min(\bar\tau,\tau')/\tau'$, 
where $\bar \tau$ is the largest number such that the restriction of $\sigma$ to $[0,T]$ has  dwell-time $\bar \tau$. 
In particular,  $\xi$ is always positive and converges to $1$ as $\tau'\searrow \tau$. 
Note that $\sigma$ and $\sigma'$ are equal except on a set of measure upper bounded by $C(T,\bar \tau)(\tau'-\tau)$, where $C(T,\bar \tau)$ denotes some positive constant only depending on $T$ and $\bar \tau$. 
As a consequence, for any $u\in L_2(0,T)$,
%\begin{equation}\label{m-eq:tau}
\[
\lim_{\tau'\searrow \tau} \Vert y_{u,\sigma'}\Vert_{2,T}=\Vert y_{u,\sigma}\Vert_{2,T}.
\]
%\end{equation}
One immediately deduces property {\bf (ii)}  from the definition of $\gamma_2(\tau,T)$ and the monotonicity of $\tau'\mapsto \gamma_2(\tau',T)$.

We can then conclude the proof of the proposition as follows: given $\tau\geq 0$ and $\kappa<\gamma_2(\tau)$, let $T>0$ be such that $\kappa<\gamma_2(\tau,T)$ ($T$ exists by property {\bf (i)}). Take then a 
right-neighbourhood $[\tau,\tau+\eta)$ of $\tau$ such that $\kappa<\gamma_2(\tau',T)$ for every $\tau'\in[\tau,\tau+\eta)$ (the existence of $\eta>0$ follows by property {\bf (ii)}). We deduce from  the monotonicity of each function $T\mapsto \gamma_2(\tau',T)$ that $\kappa<\gamma_2(\tau')$ for every $\tau'\in[\tau,\tau+\eta)$. 
Since $\kappa$ is arbitrary, we conclude that $\liminf_{\tau'\searrow \tau}\gamma_2(\tau)\ge \gamma_2(\tau)$ and the right-continuity of $\gamma_2$ is proved thanks to its monotonicity. 
\end{IEEEproof}

An immediate consequence of the monotonicity of $\tau\mapsto\gamma_2(\tau)$ and its right-continuity at $\tau=0$ is the following result. 

\begin{corollary}\label{AQ1}
The function $\gamma_2:(0,\infty)\to [0,\infty]$ is bounded if and only if $\gamma_2(0)<\infty$. 
\end{corollary}
Notice that $\gamma_2(0)$ is finite if (and only if in the uniformly
observable case)
the generalized spectral radius of a corresponding minimal realization
is smaller than one (see Theorem~\ref{th1}).

\begin{remark}\label{r-general-r-continuity}
One can reason similarly for other classes of switching laws (in particular those introduced Section~\ref{f:conca}), getting similar results to Proposition~\ref{m-prop-cont} and Corollary~\ref{AQ1}. More precisely, if we deal with a one-parameter family $[0,\bar\alpha)\ni \alpha\mapsto \mathcal{S}^\alpha(\mathscr{M})$ of classes of switching laws,
then the corresponding parameter-dependent $L_2$-gain $\gamma_2(\alpha)$ is right-continuous and non-increasing with respect to $\alpha$ provided that 
$\mathcal{S}^\alpha(\mathscr{M})\supset \mathcal{S}^{\alpha'}(\mathscr{M})$ for $\alpha<\alpha'$ and that, for $T>0$, $\alpha\in[0,\bar\alpha)$, $u\in L_2$, $\sigma\in \mathcal{S}^\alpha(\mathscr{M})$, and every sequence $\alpha_n\searrow \alpha$, there exist a sequence $(\sigma_n)_n$ with $\sigma_n\in \mathcal{S}^{\alpha_n}(\mathscr{M})$ for every $n\in\N$ and $\lim_{n\to\infty}\|y_{\sigma_n,u}\|_{2,T}=\|y_{\sigma,u}\|_{2,T}$. 
\end{remark}

 \section{Conclusion}

 In this paper, we address an open problem proposed by J.P.~Hespanha in \cite{Hespanha-unsolved}, which consists of three questions about the dependence of the $L_2$-gain $\gamma_2(\tau)$ of a switched linear control systems with respect to the dwell-time $\tau$ of its switching laws. 
 We  provide some partial answers to these questions. In particular, we prove that 
 the gain function $\tau\mapsto \gamma_2(\tau)$ is right-continuous. Further results are obtained under the assumption of uniform observability (see Definition~\ref{def:UOC}) which is equivalent, in the case of finitely many modes, to the observability of each mode of a minimal realization of the original system.
 Under such an assumption, one has that $(a)$ the infimum $\tau_{\min}$ of the dwell-times $\tau$ such that $\gamma_2(\tau)$ is finite
coincides with the largest dwell-time $\tau$ for which the 
generalized spectral radius $\rho(\tau)$  of a corresponding minimal realization is equal to one;
$(b)$ if the gain function $\gamma_2$ is bounded then 
the 
generalized spectral radius of a minimal realization corresponding to arbitrary switching is proved to be smaller than one (the converse holding true even without assuming uniform observability). 

In order to deduce these results, we are led to building an abstract framework, allowing one 
to address the same $L_2$-gain issues 
for many classes of parameter-dependent switching laws. 
The main difficulty arises from the fact that most of such classes (e.g., those of switching laws with positive dwell-time) are not closed under concatenation. To overcome this obstacle, we introduce the concepts of quasi-Barabanov semigroup and quasi-extremal trajectory.

A complete answer to the questions asked by J.~P.~Hespanha remains to be 
achieved: the continuity of the gain function $\tau\mapsto \gamma_2(\tau)$ is still open, as well as a new characterization of $\tau_{\min}$ when the uniform observability assumption does not hold. %Indeed, the above characterization of $\tau_{\min}$ may fail to be true if the uniform observability assumption does not hold. 
By the results of the paper, one knows that $\tau_{\min}$ belongs to the closed interval $I=\{\tau\geq 0\mid \rho(\tau)=1\}$. Notice, however, that $I$ may be nontrivial, as it is the case for the system considered in Example~\ref{prop-partial},
%Proposition~\ref{prop-partial}, 
for which $\tau_{\min}$ coincides with the left endpoint of $I$, contrarily to the uniformly observable case, where it is located at its right endpoint. 
It would be interesting to understand the exact location of $\tau_{\min}$ within $I$ in the general not uniformly observable case.

Another challenging open problem consists in extending the results of \cite{Margaliot-Hespanha,GeromelColaneri}, which provide an algorithmic approach based on optimal control for the computation of the $L_2$-gain, in the abstract framework introduced here.

\bibliographystyle{IEEEtran}
\bibliography{biblio-switch}

\end{document}